\documentclass[11pt]{article}
\usepackage{epsfig}
\usepackage{amssymb,amsmath}
\usepackage{epigraph}

\usepackage{graphicx}				
\usepackage{amssymb}

\newcommand{\hypgeo}[2]{%
  \operatorname{%
    {\vphantom{\mathnormal{F}}}_{#1}%
    \kern-\scriptspace
    \mathnormal{F}_{#2}%
  }%
}

\newcommand{\beq}{\begin{equation}}
\newcommand{\eeq}{\end{equation}}
\newcommand{\bea}{\begin{eqnarray}}
\newcommand{\eea}{\end{eqnarray}}

\def\boldphi{\text{\boldmath$\phi$}}

\def\boldcalD{\text{\boldmath$\mathcal{D}$}}

\begin{document}



\begin{center}

{\LARGE
SCOP: Schr{\"o}dinger Control Optimal Planning 
\vskip0.5cm
for Goal-Based Wealth Management

}

\vskip1.0cm
{\Large Igor Halperin} \\
\vskip0.5cm
ighalp@gmail.com\footnote{Opinions presented in this paper are author's only, and not necessarily  of his employer. The standard disclaimer applies.}\\
\vskip0.5cm
\today \\

\vskip1.0cm
{\Large Abstract:\\}
\end{center}
\parbox[t]{\textwidth}{
We consider the problem of optimization of contributions of a financial planner such as a working individual towards a financial goal such as retirement. The objective of the planner is to find an optimal and feasible schedule of periodic installments to an investment  portfolio set up towards the goal.  Because portfolio returns are random, the practical version of the problem amounts to finding an optimal contribution scheme such that the goal is satisfied at a given confidence level. This paper suggests a semi-analytical approach to a continuous-time version of this problem based on a controlled backward Kolmogorov equation (BKE) which describes the tail probability of the terminal wealth given a contribution policy. 
The controlled BKE is solved semi-analytically by reducing it to a controlled Schr{\"o}dinger equation and solving the latter using an algebraic method.
Numerically, our approach amounts to finding semi-analytical solutions simultaneously for all values of control parameters on a small grid, and then using the standard two-dimensional spline interpolation to simultaneously represent all satisficing solutions of the original plan optimization problem. 
Rather than being a point in the space of control variables, satisficing solutions form continuous contour lines (efficient frontiers) in this space.
  }
 \newcounter{helpfootnote}
\setcounter{helpfootnote}{\thefootnote} 
\renewcommand{\thefootnote}{\fnsymbol{footnote}}
\setcounter{footnote}{0}
\footnotetext{
I thank Ernest Baver, Shaina Race Bennett, Eric Berger, Lisa Huang, Andrey Itkin, Sergey Malinin, Mirco Miletari and Jack Sarkissian for helpful discussions.
}     

 \renewcommand{\thefootnote}{\arabic{footnote}}
\setcounter{footnote}{\thehelpfootnote} 

\newpage
 
\section{Introduction}

We consider the problem of optimization of contributions of a financial planner such as a working individual towards a financial goal such as retirement. The objective of the planner is to find an optimal and feasible schedule of periodic installments to an investment  portfolio set up towards the goal. 
Because portfolio returns are random, in practice this problem amounts to finding an optimal contribution scheme such that the goal is satisfied at a given confidence level. This problem is similar to  the celebrated Merton (1971) optimal consumption problem \cite{Merton_1971}, however it differs from the latter in at least three important aspects. 

First, we do not address a portfolio optimization part of this problem, instead assuming 
that an investment portfolio set towards the goal has a fixed equity to cash ratio which is maintained by a financial fiduciary. This implies that a part of a regular cash contribution to the portfolio from the planner that is invested in equity is fixed by the fiduciary and is not subject to optimization.\footnote{Our approach can also be extended to incorporate portfolios with a decreasing share of equity to mimic glide path profiles popular in the industry.} 

Second, we do not maximize any utility of consumption, and respectively do not consider any consumption process either. Instead, we work directly with contributions to the financial plan as decision variables constrained by a net income, and our objective is to exceed a certain target wealth level at a terminal time $ T $ at a given confidence level $ \alpha $. In this sense, our problem formulation fits the setting of goal-based wealth management, see e.g. Das {\it et. al.} (2018) \cite{Das_2018} and references therein.      

The problem as formulated belongs in the class of optimal stochastic control problems, and thus could in principle be solved using conventional tools such as 
dynamic programming (DP) or reinforcement learning (RL).\footnote{See e.g. Dixon, Halperin and Bilokon (2020) \cite{MLF} for a review with financial applications including applications to goal-based financial planning.} However, our problem is not necessarily a problem of {\it adaptive} (closed-loop) control type that are usually addressed with DP and RL methods. Retirement planning usually involves setting up a plan upfront, which is then revisited, if necessary, on the annual basis. A typical plan usually involves a time lime of annual installments, where annual contributions planned for the future are scheduled to increase to adjust for inflation, as well as possibly incorporate 
benefits from expected future increases in wages or other investable income. Such settings are typical for {\it open-loop} control problems, where an agent decides on the policy at the start, and then keeps it fixed for a certain fixed number of steps, without adapting it at each steps as is done in DP or RL.

Another important difference of our problem from the setting of DP or RL is that the latter methods solve a maximization problem for a value (or action-value) function. If the value function has a global maximum as a function of a deterministic policy, then there should be a policy - a point in a multi-dimensional action space whose position depends on the state - that attains this maximum.
In contrast, a proper framework for decision-making in our setting is the framework of {\it satisficing decision making} of Herbert Simon (1955, 1979( \cite{Simon_1955, Simon_1979}. 
With this framework, we do not maximize an objective function. Instead, we only want to bring it to a certain {\it goal level} that we define beforehand. This involves two steps: {\it search} (i.e. exploring the objective function for different values of input controls), and {\it satisficing} (i.e. picking a solution that satisfies the goal criterion) \cite{Simon_1979}. Clearly, with this setting, feasible problems should in general produce continuous lines in the action space, rather than isolated points, as optimal solutions. 
 
This paper suggests a semi-analytical solution to a continuous-time version of the problem formulated above. 
We use a probabilistic approach that amounts to analysis of two related equations. The first one is Fokker-Planck equation (FPE) whose solution is the probability density of the terminal wealth at the at time $ T $, for a given contribution plan. The second equation is the backward Kolmogorov equation (BKE) which gives the probability of the terminal wealth to fall below a certain target wealth at a given confidence level $ \alpha $.  
We solve the resulting controlled FPE and BKE equations semi-analytically by reducing them to controlled Schr{\"o}dinger equations\footnote{The Schr{\"o}dinger equations obtained in our analysis are controlled in the sense that that their parameters and initial positions of a 'particle' representing the planner's wealth are controlled by the user.}, and then solving the latter equations in a semi-closed form using an algebraic method. These semi-analytical solutions are computed {\it simultaneously} for all values of two control parameters describing annual contributions and their growth rate, on a small 2D grid of their discretized admissible values. Once the values of the tail BKE probability are computed for all nodes on the grid, its values for arbitrary inputs are cheaply obtained using the standard 2D cubic spline interpolation. This completes the {\it search} step of satisficing decision-making in our problem. 

Once the {\it search} stage is completed, the {\it satisficing} part is rather simple, as all it takes is to draw constant-probability lines at the confidence levels $ \alpha $ on the surface of the BKE tail probability function viewed as a function of initial controls. This task can be done using the standard 2D spline interpolation plus a root-finding algorithm, or equivalently using a contour-plotting algorithm that already combines the first two algorithms.   

Therefore, in our approach, dubbed the Schr{\"o}dinger Control Optimal Planning, or SCOP for short, the initial problem of contribution optimization amounts to the standard 2D spline interpolation plus a root-solving (or equivalently contour-plotting) algorithm, whose end result gives a visual representation of {\it all} possible satisficing plans at different confidence levels.
Such visual representation is conceptualized via the notion of {\it efficient frontiers}: constant-level lines on the 2D surface of the BKE probability as the function of controls. This gives the planner the ability to perform a real time policy optimization for a given goal, or scenario analyses for different settings of the planning problem.

The paper is organized as follows. In Sect.~\ref{sect_OC}, we present the theoretical framework for our problems.
Sect.~\ref{sect_plan_optimization} describes our method for retirement plan optimization.
Finally, a summary and an outlook for future research are presented in Sect.~\ref{sect_Summary}.  All technical details and derivations are left to Appendices A, B, and C. The appendices have a hierarchical structure, i.e. Appendices B and C clarify certain technical details for topics discussed in Appendix A.

\section{Optimal contribution problem}
\label{sect_OC}

\subsection{The investment portfolio dynamics}
\label{sect_portfolio_dynamics}

We consider a financial planner such as a retirement planner with a planning horizon $ T $ who invests in a portfolio created towards a goal, expressed as a desired terminal wealth at time $ T $. The investment portfolio has a fixed ratio $ \omega $ of equity to the total portfolio value (e.g. $ \omega = 0.9 $), which is assumed to be maintained by a fiduciary that manages the portfolio on client's behalf. 

Let $ c_t \Delta t $ be the after-tax money flow from the planner to the portfolio in the time step $ [t, t + \Delta t] $. This investment of cash from the planner forces the fiduciary to buy the amount of equity equal to $ \omega c_t \Delta t $, which also incurs proportional transaction cost $ \nu \omega | c_t | \Delta t $, where $ \nu $ is a parameter. As we have $ c_t \geq 0 $ in our problem, in what follows we omit the absolute value expression in the transaction cost formula. 

The investment portfolio is composed of cash $ b_t $ and equity $ s_t $. We assume that the cash infusion to the portfolio happens at time $ t_{+} $ immediately after $ t $, followed by an immediate purchase of equity by the fiduciary. The cash and equity positions at time $ t_{+} $ are then obtained as follows:
\bea
\label{cash_equity_t+}
b_{t}^{+} &\hspace{-0.2cm} = \hspace{-0.2cm} & b_t + \left((1-\omega) c_t  -\nu \omega c_t \right) \Delta t \nonumber \\
s_t^{+} &\hspace{-0.2cm} = \hspace{-0.2cm} & s_t + \omega c_t \Delta t
\eea
This produces the following update rule for the total portfolio value $ \Pi_t := b_t + s_t $:
\beq
\label{Pi_update}
\Pi_{t}^{+} = \Pi_t + \left(c_t - \nu \omega c_t \right) \Delta t
\eeq
We assume that equity has random returns
\beq
\label{equity_ret}
r_e(t) = \hat{r}_e + \frac{\hat{\sigma}_e}{\sqrt{\Delta t}} \xi_t, \; \; \; \xi_t \sim \mathcal{N}(0,1)
\eeq
with a mean value $ \hat{r}_e $, while the cash component has a deterministic return $ r_f $. Using the fact that $ s_t = \omega V_t $ and taking the continuous-time limit $ \Delta t = d t \rightarrow 0 $, we obtain a stochastic differential equation (SDE) for the portfolio value $ V_t $:
\beq
\label{SDE_port_Pi}
d \Pi_t = \left( \bar{r} \Pi_t + u_t \right) dt + \sigma \Pi_t d W_t
\eeq
here $ W_t $ is a standard Brownian motion, and the decision variable $ u_t $ and parameters $ \bar{r}, \sigma $ are defined as follows:
\beq
\label{bar_r_u_t}
 u_t := (1 - \nu \omega) c_t , \; \; \; \bar{r} := r_f + \omega (\hat{r}_e - r_f), \; \; \; \sigma := \omega \hat{\sigma}_e 
 \eeq      
In this paper, we consider contribution policies parameterized in terms of an initial contribution rate $ u_0 $ and a growth rate $ \xi $:
\beq
\label{u_t_class}
u_t = u_0 e^{\xi t}
\eeq
so that $ u_0 $ and $ \xi $ are now the constrained decision variables for the planner. The growth rate $ \xi $ can be tuned to compensate for inflation or 
increase future contributions, assuming that the net investable income will grow in the future. 

\subsection{Stochastic Verhulst equation and Langevin dynamics}
\label{sect_Langevin}

Note that Eq.(\ref{SDE_port_Pi}) is a controlled SDE where control variables are $ u_0 $ and $ \xi $. To simplify these dynamics, we introduce a new dimensionless variable $ v_t $ defined as follows:
\beq
\label{y_Morse_var}
v_t = \frac{2}{\sigma^2} \frac{u_t}{\Pi_t} = \frac{2}{\sigma^2}  \frac{u_0 e^{\xi t}}{\Pi_t}, \; \; \; 
\eeq 
Note that the new variable $ v_t $ blends together the state variable $ \Pi_t $ and action variable $ u_t $. 
Using It{\^o}'s lemma, we obtain the SDE for $ v_t$:
\beq
\label{SDE_Morse}
d v_t = v_t \left( \xi - \bar{r} + \sigma^2 - \frac{\sigma^2}{2} v_t \right) dt +  \sigma v_t d W_t
\eeq
The new SDE for $ s_t $ does not depend on $ u_0 $, which is now embedded into the problem via the initial value $ v_0 $. The only remaining control parameter explicitly appearing in (\ref{SDE_Morse}) is $ \xi $.

The SDE (\ref{SDE_Morse}) is known in the literature as the stochastic Verhulst equation, see e.g. Giet {\it et. al.} (2015) \cite{Giet_2015}. Recently, the SDE (\ref{SDE_Morse}) with time-dependent coefficients was proposed by Itkin, Lipton and Muravey (2020) \cite{Itkin} as an attractive model of short interest rates that addresses some deficiencies of the classical Black-Karasinski model.
Furthermore, they developed a method of solving the BKE for the SDE (\ref{SDE_Morse}) with general time-dependent coefficients using an original integral transform method \cite{Itkin_book}. The method of Itkin, Lipton and Muravey (2020) \cite{Itkin} could be considered as an alternative approach to the semi-analytical method developed in this paper, see also the summary section for further comments.\footnote{The author would like to thank Andrey Itkin for conversations about 
the stochastic Verhulst model and mathematical connections between this work and \cite{Itkin}.}  
 
To transform the multiplicative noise in the resulting SDE into an additive noise, we apply one more transformation to a new variable $ x_t $:
\beq
\label{x_var}
x_t = - \log v_t 
\eeq
so that $ v_t = e^{-x_t} $. 
The SDE for $ x_t $ reads
\beq
\label{SDE_x}
d x_t = - \left( \eta - \frac{\sigma^2}{2} e^{-x} \right) dt + \sigma d W_t, \; \; \; \eta := \xi - \bar{r} + \frac{\sigma^2}{2}
\eeq
The new parameter $ \eta $ introduced in the SDE (\ref{SDE_x}) can be used as a new control variable in liew of $ \xi $, with constraints inherited from constraints on $ \xi $. 

Now, the SDE (\ref{SDE_x}) can be viewed as an (overdamped) Langevin equation 
\beq
\label{Langevin}
d x_t =  - \frac{\partial V}{\partial x_t} dt + \sigma d W_t
\eeq
where the potential function $ V = V(x) $ is defined as follows:
\beq
\label{V_x}
V(x) = 
 \eta x + \frac{\sigma^2}{2} e^{-x}
\eeq
so that $ \eta $ becomes a control parameter for the Langevin potential $ V(x) $.

Consider now the conditional tail probability $ p(w,t) := P^{(\pi)}[\Pi_T < \hat{\Pi}| \Pi_t = w] $,  where $ \Pi_T $ stands for the terminal wealth at time $ T $, 
$ \hat{\Pi} $ is a target wealth, and the notation $ P^{(\pi)} \left[ \cdots \right] $ is introduced to emphasize the dependence of this conditional probability on the policy $ \pi := (u_0, \xi) $.
 Our problem amounts to the task of choosing optimal values of $ u_0 $ and $ \xi $ such that the conditional tail probability $ p(v,0) $ as seen now at time $ t = 0 $ does not exceed a certain target probability level $ \alpha $ such as $ \alpha = 5\% $. Now, when conditioning on the policy  $ \pi := (u_0, \xi) $, we can equivalently compute the tail probability $ p(v,t) $ in terms of variables $ v_t $ and $ x_t $ introduced above. Denoting corresponding target variables computed according to 
 Eqs.(\ref{y_Morse_var}) and (\ref{x_var}) as $ \hat{v} $ and $ \hat{x}$, we can write down the following chain of equations:
 \beq
 \label{chain_of_tail_probs_vyz}
 P^{(\pi)} \left[ \Pi_T \leq \hat{\Pi}| \Pi_t = w \right] = P^{(\pi)} \left[ v_T \geq \hat{v} | v_t = v \right] = P^{(\pi)} \left[ x_T \leq \hat{x} | x_t = x \right]
 \eeq   
This shows that the tail probability $ P^{(\pi)} \left[ \Pi_T \leq \hat{\Pi}| \Pi_t = w \right] $ of the original problem can be equivalently computed using the 
variable $ x_t $. As will be shown next, the advantage of doing this is that using the variable $ x_t $, our problem has a semi-analytical solution.    


\subsection{From the Fokker-Planck-Kolmogorov control to the Schr{\"o}dinger control }
\label{sect_FP}

The Langevin dynamics expressed by Eq.(\ref{Langevin}) have an equivalent probabilistic formulation in terms of the corresponding Fokker-Planck equation (FPE) 
for the probability density $ f(x,t) $:
\beq
\label{FPE}
\frac{\partial f(x,t)}{\partial t} =  \frac{\partial}{\partial x} \left[ \frac{\partial V}{\partial x} f(x,t) \right] + \frac{\sigma^2}{2} \frac{\partial^2 f(x,t)}{\partial x^2}
\eeq
The probability density $ f(x,t) $ that solves the FPE equation is a distribution of $ x_t $ at time $ t $ given controls $ (u_0, \eta) $. To ease the notation, the dependence on controls is suppressed in all relations in this section. 

The second probability function we consider is the tail probability $ p(x,t) := P^{(\pi)} \left[ x_T \leq \hat{x} | x_t = x \right] $. It satisfies the backward Kolmogorov equation (BKE):
\beq
\label{BKE}
\frac{\partial p}{\partial t} - \frac{\partial V}{\partial x} \frac{\partial p}{\partial x} + \frac{\sigma^2}{2} \frac{\partial^2 p}{\partial x^2} = 0
\eeq 
with the terminal condition
\beq
\label{terminal_Kolm}
p(x,T) = \theta \left(\hat{x} - x \right), \; \; \;  \hat{x} := 
\log \frac{ \hbar \hat{\Pi} }{ 2 u_0 e^{\xi T}} 
\eeq 
Let us look for solutions of the FPE (\ref{FPE}) and BKE (\ref{BKE}) in the following form:
\bea
\label{Schrodinger_change}
f(x,t) &\hspace{-0.2cm} = \hspace{-0.2cm} & e^{ - \frac{V(x)}{\hbar} } \Psi(x,t) , \; \; \; \hbar := \sigma^2  \\ 
p(x,t) &\hspace{-0.2cm} = \hspace{-0.2cm} & e^{\frac{V(x)}{\hbar}} \bar{\Psi}(x,\tau), \; \; \; \tau := T - t\nonumber
\eea
Here we introduced a new parameter $ \hbar = \sigma^2 $ to make formulae to follow look more familiar. Also note the time reversal $ t \rightarrow \tau = T - t $ in the second relation. 

Substituting these expressions back into Eqs.(\ref{FPE}) and (\ref{BKE}), we obtain a pair of PDEs  for $ \Psi(z,t) $ and $ \bar{\Psi}(z,\tau) $:
\beq
\label{SE}
- \hbar \frac{\partial \Psi(x,t)}{\partial t} = \left[ - \frac{\hbar^2}{2} \frac{\partial^2}{\partial x^2} + U(x) \right] \Psi(x,t)
\eeq
\beq
\label{SE_back}
 -\hbar \frac{\partial \bar{\Psi}(x,\tau)}{\partial \tau} = \left[ - \frac{\hbar^2}{2} \frac{\partial^2}{\partial x^2} + U(x) \right] \bar{\Psi}(x,\tau)
\eeq
where 
\beq
\label{U_x}
U(x) = \frac{1}{2} \left( \frac{\partial V}{\partial x} \right)^2  - \frac{\hbar}{2} \frac{\partial^2 V}{\partial x^2}
\eeq
 Written in this form, Eqs.(\ref{SE}) and (\ref{SE_back}) should look very familiar to the reader acquainted with quantum mechanics.
Eq.(\ref{SE}), formulated in time $ t $, is the Schr{\"o}dinger equation in Euclidean (imaginary) time with the 'Planck constant' $ \hbar = \sigma^2 $ and potential
 $ U(x) $.   On the other hand, Eq.((\ref{SE_back}) is the Schr{\"o}dinger equation in backward time $ \tau = T - t $. 
 In what follows, we will refer to Eqs.(\ref{SE}) and (\ref{SE_back}) as, respectively, the Fokker-Plack Schr{\"o}dinger equation (FP-SE), and the backward Kolmogorov Schr{\"o}dinger equation (BK-SE).
 
 Note that while the FP-SE (\ref{SE}) and BK-SE (\ref{SE_back}) have identical forms, they satisfy different initial conditions:
 \beq
 \label{init_cond_Psi}
 \Psi(x,0) = e^{\frac{V(x)}{\hbar}} \delta \left(x - x_0 \right), \; \; \; \;
 \bar{\Psi}(x,0) = e^{- \frac{V(x)}{\hbar}} \theta\left( \hat{x} - x  \right)
 \eeq
 This implies that $ e^{- \frac{V(x)}{\hbar}} \Psi(x,t) $ is the Green's function for  Eq.(\ref{SE_back}), and therefore, we can relate the two functions as follows:
 \beq
 \label{Psi_bar_Psi}
 \bar{\Psi}(x,T) = \int_{- \infty}^{\infty} dx' e^{- \frac{V(x')}{\hbar}}  \bar{\Psi}(x',0) \Psi(x,T| x',0)
 \eeq
 Using Eq.(\ref{Schrodinger_change}), this relation is seen to be equivalent to the conventional Feynman-Kac expression for the tail probability $ p(x,0) $:
 \beq
 \label{tail_prob_BKE}
 p(x,0) =   e^{\frac{2 V(x)}{\hbar}} \int_{- \infty}^{\infty} dx' e^{-  \frac{2 V(x')}{\hbar}} \theta(\hat{x} - x') f(x,T| x',0)
 \eeq
 
 While the Feynman-Kac formula (\ref{tail_prob_BKE}) gives one representation of a solution to our problem, an alternative approach is to directly solve the SE (\ref{SE_back}) with the initial condition stated in (\ref{init_cond_Psi}). As Eqs.(\ref{SE}) and (\ref{SE_back}) have the same Hamiltonian, their solutions can be found by expanding in the same set of basis functions. Therefore, it might be more convenient to directly solve the the BK-SE (\ref{SE_back}) because it does not involve an extra integration step of the Feynman-Kac representation (\ref{tail_prob_BKE}), at least formally.\footnote{We say here 'at least formally' because the Feynman-Kac integration step in Eq.(\ref{tail_prob_BKE}), instead of being performed at the final step as in (\ref{tail_prob_BKE}), in fact has an exact parallel step within an eigenvalue decomposition for (\ref{SE_back}): weights for the eigenvalue decomposition for the backward and forward problems are related by the same integral
 relation (\ref{tail_prob_BKE}) as the probabilities themselves. For more details, see Appendix A.
 }  
 
 \subsection{The Morse potential}
 \label{sect_Morse_potential}
 
 The Schr{\"o}dinger potential (\ref{U_x}) for our problem is obtained using Eq.(\ref{V_x}):
 \beq
 \label{U_Morse_x}
 U(x) = 
 D_0 \left( 1 - \frac{1}{g} e^{-x} \right)^2 + U_0  
 \eeq
where we defined the following parameters:
\beq
\label{params}
g := 1 + \frac{2 \eta}{\hbar} = 2 \left( 1 + \frac{\xi - \bar{r}}{\hbar} \right), \; \; \;
D_0 := \frac{\hbar^2}{8}g^2, \; \; \;   
U_0 :=  - \frac{\hbar^2}{8}(2g-1)
 \eeq
The control variable $ \xi $ of the original problem is now embedded into the 'coupling constant' parameter $ g  $ of the corresponding quantum mechanical problem
(\ref{SE}) or (\ref{SE_back}). The original control problem therefore becomes a problem of {\it quantum control} for the Euclidean quantum mechanics. In this problem, we control the coupling constant $ g $ of the quantum mechanical potential $U(x)$ and the initial position of a quantum mechanical 'particle' representing the wealth of the planner in such a way that the terminal WF $ \bar{\Psi}(x,T) $ is driven to a desired form. 

 In its turn, the coupling constant $ g $ determines the degree of non-linearity of potential $ U(x) $. If $ g > 0 $,  the potential has a minimum at $ x =  - \log g $. 
In the limit $ g \rightarrow 0 $, the minimum is pushed to infinity, and we obtain $ U(x) = (\hbar^2/8) \left( 1 + e^{-2x} \right) $, i.e. the potential becomes a monotonically decreasing function. In the opposite limit $ g \rightarrow \infty $, the potential $ U(x) $ becomes a constant $ U(x) = D_0 + U_0 $. For $ g < 0 $, the potential $ U(z) $ has no minima or maxima, and becomes a monotonically decreasing function. 

\begin{figure}[ht]
\begin{center}
\includegraphics[
width=95mm,
height=60mm]{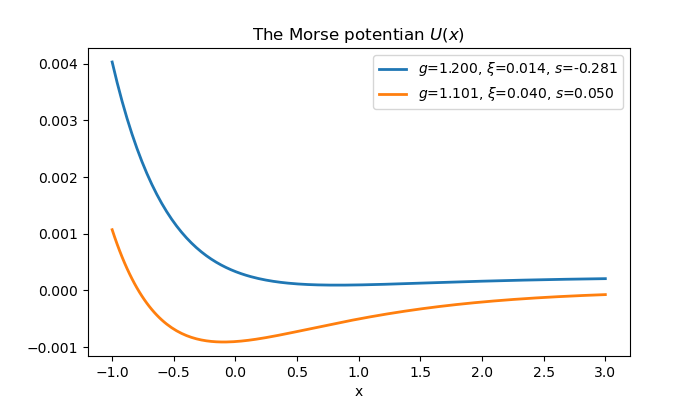}
\caption{Examples of the Morse potential for different values of parameter $ g $ driven by different values of control parameter $ \xi $. Both potentials have a minimum at $ x_{\star} = - \log(g) $ where $ g = 2s+1 $. For small values of $ \xi $, the potential well becomes very shallow. If $ \xi $ becomes smaller than 
$ \xi_{low} := \bar{r} - \sigma^2 $, the minimum disappears, and the potential becomes monotonically decreasing. For larger values of $ \xi $, the minimum becomes more pronounced.}  
\label{fig_Morse_potential}
\end{center}
\end{figure}

The potential $ U(x) $ can be written more compactly once we differentiate between two possible regimes $ g > 0 $ and $ g < 0 $ of the model.
Introducing parameter $ x_{\star} := - \log | g| $,  we can write the potential $ U(x) $ in the following form: 
\beq
 \label{U_Morse}
 U(x) = 
 \left\{ \begin{array}{l}
D_0 \left(1 - e^{-(x-x_{\star})} \right)^2 + U_0 , \; \; \; \text{if} \; \; \; g  > 0 \\       
D_0 \left(1 + e^{-(x - x_{\star} )} \right)^2 + U_0, \; \; \; \text{if} \; \; \; g   <  0  \\  
 \end{array} \right.
 \eeq
 Focusing on the first expression here obtained for $ g > 0 $, this is the celebrated Morse potential known from quantum mechanics, see Landau (1980) \cite{Landau}, see Fig.\ref{fig_Morse_potential}. As this potential is very tractable, the reduction of our initial financial problem to Schr{\"o}dinger equations (\ref{SE}), (\ref{SE_back}) with the Morse potential $ U(x) $ is very helpful. More specifically, it enables semi-analytical solutions for both 
 Schr{\"o}dinger equations (\ref{SE}) and (\ref{SE_back}), and hence for both the FPE (\ref{FPE}) and BKE (\ref{BKE}) by virtue of Eq.(\ref{Schrodinger_change}).\footnote{The link between the FPE equation (\ref{FPE}) and the  Schr{\"o}dinger equation (\ref{SE}) via Eq.(\ref{Schrodinger_change}) is well known in the literature, see e.g. Van Kampen (1982) \cite{vanKampen}.  The Euclidean-time SE (\ref{SE}) is obtained from the conventional quantum mechanical SE by the so-called Wick rotation $ t \rightarrow - i t $ and then taking $ t $ to lie in an interval $ t \in [0,T] $.  
 Similarly, the second Euclidean-time equation (\ref{SE_back}) can be obtained by using instead a different substitution $ t \rightarrow -i (T-t) $ and then taking the variable $ \tau := T-t $ to lie in the interval $ \tau \in [0, T] $. The only difference from the first Euclidean SE is that now an initial condition for $ \bar{\Psi}(x,t)$ is determined by a terminal condition for the BKE probability $ p(x,T)$. As in physics we deal with forward problems, such alternative Wick rotations are usually not considered there, however it is of interest in our backward control problem. To our knowledge, the second transformation in Eq.(\ref{Schrodinger_change}) for the BKE equation was not previously considered in the stochastic control literature.} We 
 only give here the final expressions for 
the FPE and BKE solutions $ f(x,t) $ and $ p(x,t) $. Details are left to Appendix A which is written in a reasonably self-contained way to provide a brief overview of these methods for the reader without a background in theoretical physics. 

 The solution presented in Appendix A shows that both functions $ \Psi(x,t)$ and $ \bar{\Psi}(x,t) $ depend on $ x $ only through the so-called Morse 
 variable $ y := e^{-x} $, so that we can equivalently write them as $ \Psi(y,t), \, \bar{\Psi}(y,t) $. Recalling Eq.(\ref{x_var}), the time-$ t $ Morse variable value $ y_t $ can be expressed in terms of the original model inputs:
 \beq
 \label{Morse_var}
 y_t = e^{-x_t} = \frac{ 2 u_0 e^{\xi t}}{ \hbar \Pi_t}
 \eeq 
 If we use this relation for $ t = 0 $, it means that the initial value $ y_0 $ at time zero is simply the initial contribution $ u_0 $ scaled by the portfolio value $ \Pi_0 $ and portfolio variance $ \hbar = \sigma^2 $.
 In other words, $ y_0 $ can be viewed as a correct dimensionless version of the initial control $ u_0 $ (which has the dimension of dollars per year) which enables a meaningful comparison of different plans and different market conditions. For this reason, in most of our numerical examples to follow we will use the initial value  $ y_0 $ of the Morse variable as a dimensionless control parameter in lieu of the initial control parameter $ u_0 $. 
 
 As $ \Psi(x,t), \, \bar{\Psi}(x,t) $ are functions of the Morse variable $ y $,  it is convenient to reformulate Eqs.(\ref{Schrodinger_change}) in terms of the Morse variable $ y = e^{-x} $. For the exponential factor $ e^{-V(x)/\hbar} $, it gives
 \beq
 \label{exp_V_in_y}
 e^{-V(y)/\hbar} = \left. e^{- V(x)/\hbar} \right|_{x=-\log y} = y^{\eta/\hbar} e^{-y/2}
 \eeq  
Changing the variables in Eq.(\ref{Schrodinger_change}) to the Morse variable $ y = e^{-x} $, we write them as follows\footnote{The additional factor 
$ y^{-1} $ in the first equation in (\ref{Schrodinger_change_y}) is the Jacobian of transformation $ x \rightarrow y $.}
\bea
\label{Schrodinger_change_y}
f(y,t) &\hspace{-0.2cm} = \hspace{-0.2cm} & y^{-1} e^{ - V(y)/\hbar } \Psi(y,t) = y^{\eta/\hbar-1} e^{-y/2} \Psi(y,t)   \\ 
p(y,t) &\hspace{-0.2cm} = \hspace{-0.2cm} & e^{ V(y)/\hbar} \bar{\Psi}(y,\tau) = y^{-\eta/\hbar} e^{y/2} \bar{\Psi}(y,\tau) \nonumber 
\eea
The explicit form of the FP-SE WF $ \Psi(y,t) $ and the BK-SE WF $ \bar{\Psi}(y,t) $ is given in Eqs.(\ref{FPE_f_final_A})  and (\ref{SE_back_res}), respectively.

\begin{figure}[ht]
\begin{center}
\hspace*{-2cm}  \includegraphics[
width=190mm,
height=85mm]{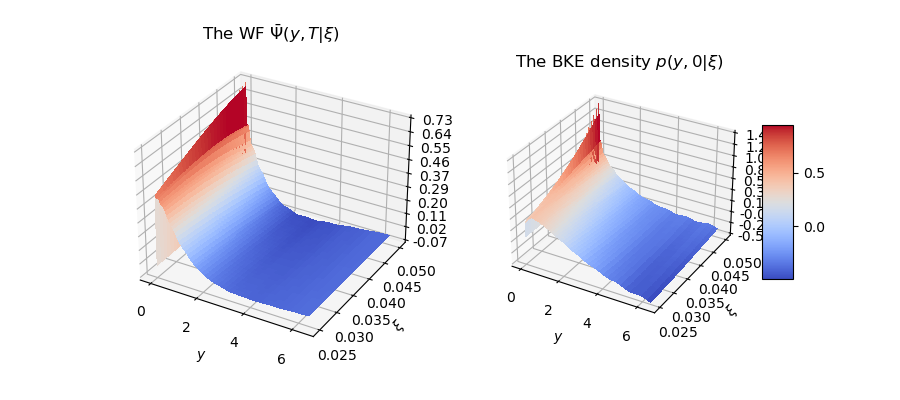}
\caption{On the left: The WF $ \bar{\Psi}(y, T) $ for $ T = 20 $, shown as a function of both $ y $ and $ \xi $ using $ N = 150 $ basis functions, for a retirement plan with the horizon $ T = 20Y $, the target wealth of \$2.5m, the initial wealth of \$0.5m, and portfolio volatility of 30\%, with the following bounds on admissible values of $ u_0 $ and $ \xi $: \$10k $\leq u_0 \leq $ \$100K, and $ 2.5\% \leq \xi \leq 5.0\% $. On the right:
the corresponding BKE tail probability $ p(y,0) $.} 
\label{fig_Psi_bar_and_BKE_3D}
\end{center}
\end{figure}

In Fig.~\ref{fig_Psi_bar_and_BKE_3D}, we show 3D views of the WF $ \bar{\Psi}(y,T)$ and the BKE density $ p(y,0) $ as a function of two variables $ y = y_0 $ and $ \xi $ for a typical retirement plan with the horizon $ T = 20Y $ for the target wealth of \$2.5m, initial wealth of \$0.5m, and portfolio volatility of 30\%, with the following bounds on admissible values of $ u_0 $ and $ \xi $:  \$10k $ \leq u_0 \leq $ \$100K, and $ 2.5\% \leq \xi \leq 5.0\% $. This parameter and constraints settings will also be used below in other examples illustrating the working of the model.  


Note that using Eqs.(\ref{Schrodinger_change_y}), we can obtain the Chapman-Kolmogorov equation in terms of $ \Psi(x,t) $ and $ \bar{\Psi}(x,t) $: 
\beq
\label{norm_psi_psibar}
\int p(x,t) f(x,t|x_0,t_0) dx = \int \Psi(x,t) \bar{\Psi}(x,t) dx = P \left[ \left. x_T \leq \hat{x} \right| x_0, t_0 \right] 
\eeq
which resembles the normalization condition $ \int \left| \Psi(x,t) \right|^2 dx = 1 $ in quantum mechanics. 

\section{Optimization of contribution policy: efficient frontiers}
\label{sect_plan_optimization}

%

Optimization of contribution policy is performed using the semi-analytical solution given by the second equation in Eqs.(\ref{Schrodinger_change_y}) which gives the tail probability of the event $ \Pi_T < \hat{\Pi} $ in terms of the tail probability of the Morse variable $ y $:
\beq
\label{reduction_for_P}
P \left[ \left.  \Pi_T \leq \hat{\Pi} \right| \Pi_t, u_t \right] = P \left[ y_T \geq \hat{y} | y_t, \xi \right] = y^{ -\eta/\hbar} e^{ y/2} \bar{\Psi}(y,t | \xi)
\eeq
Here $ \hat{y} $ is the target value of the Morse variable $ y $, which is computed according to Eq.(\ref{Morse_var}):
\beq
\label{hat_y_main}
\hat{y} = \frac{2 u_0 e^{\xi T} }{ \hbar \hat{\Pi}}
\eeq
and we wrote the WF $ \bar{\Psi}(y,t) $ as $ \bar{\Psi}(y,t | \xi) $ to emphasize its dependence on $ \xi $.

Our strategy for optimization of contribution policy is to compute Eq.(\ref{reduction_for_P}) for all possible values of control variables $u_0 $ and $ \xi $, and then find the optimal control that matches the threshold value for the BKE tail probability $ p(y,0) $. In practice, this involves computing the tail probability (\ref{reduction_for_P}) for 
all nodes of a certain two-dimensional grid for $ u_0 $ and $ \xi $, and then interpolating to other values of $ u_0 $, $ \xi $ using splines.  
As in the range of input parameter that have practical interest the model behaves in a very smooth way (see Fig.~\ref{fig_Psi_bar_and_BKE_3D}), this means that we can use a rather small grid of pairs $ (u_0, \xi) $ to perform accurate spline interpolation. 

Furthermore, the mathematical structure of the solution to our problem suggests that instead of considering the pair $ (u_0, \xi ) $ as independent controls, we can equivalently but more conveniently map them onto the pair $ (y, \xi) $.  
As suggested by Eq.(\ref{Morse_var}), $ y $ is the right dimensionless version of $ u_0 $, additionally scaled by the portfolio volatility, therefore using variable $ y $ instead of the absolute dollar amount $ u_0 $ enables comparing different plans and different market conditions.  
If the current time is set to zero and the initial wealth $ \Pi_0 $ is fixed, then values of $ u_0 $ corresponding to values of $ y = y_0 $ are obtained according to (\ref{Morse_var}):
 \beq
\label{y0_from_u0}
u_0 = \frac{\hbar}{2}   y_0 \Pi_0
\eeq
In our numerical examples presented below, we use a non-uniform grid\footnote{We use a log-spaced grid to put more points in the region of small $ y $ to better capture strong variations of FPE and BKE densities for small values of $ y $.}  of values of $ y $ of size 100. In addition, we use 20 grid points to represent possible values of $ \xi $. The values of the BKE density for arbitrary intermediate inputs are obtained using the standard 2D cubic splines. 


As the inputs $ (\xi, y) $ (or equivalently $ (\xi,u_0) $) form a 2D plane, 
Eq.(\ref{reduction_for_P}) implies that the problem of finding the optimal contributions $ u^{\star} = (u_0^{\star}, \xi^{\star} ) $ for a plan with the goal wealth $ 
\hat{\Pi} $ that we want to exceed with probability equal to $ 1- \alpha $ 
amounts to solving the equation
\beq
\label{min_algo}
P \left[ y_T > \hat{y} | y_t, \xi_t, t \right] =  \alpha
\eeq
in this 2D $ (\xi y) $-plane. 

Note that we call here an optimal solution $ u^{\star} = (u_0^{\star}, \xi^{\star} ) $ is in fact a {\it satisficing} solution: we do not want to drive the tail probability
(\ref{min_algo}) all the way down to zero, we only need to drive it below the target confidence level $   \alpha $. Satisficing decision-making is different from decision-making based on maximization of a total reward as is done in DP or RL. Indeed, with DP or RL, a reasonably-behaved value function should have a global maximum, therefore there should be at least one policy that leads to such maximum value. Quite differently, with satisficing decision-making problems such 
as our Eq.(\ref{min_algo}), solutions for feasible problems form contour lines on the surface of BKE density $ p(y,0| \xi) $.    

\begin{figure}[ht]
\begin{center}
\hspace*{-2cm} \includegraphics[
width=140mm,
height=80mm]{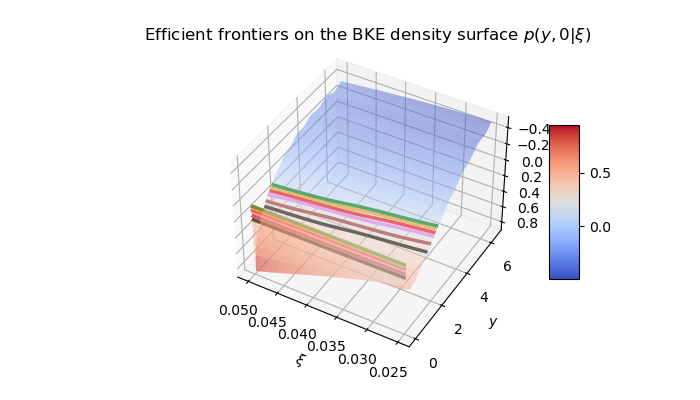}
\caption{  
3D plot of the BKE tail probability $ p(y,0| \xi) $ viewed as a function of both control variables $y $ and $ \xi $, along with efficient frontier lines corresponding to probability confidence levels $ \alpha = [3,5,7.5,10,15,20] \% $.
} 
\label{fig_BKE_frontier_3D}
\end{center}
\end{figure}

\begin{figure}[ht]
\begin{center}
\hspace*{-1cm} \includegraphics[
width=100mm,
height=65mm]{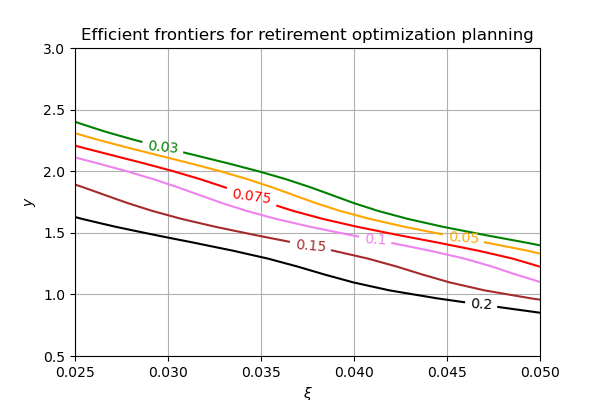}
\caption{  
The 2D view of efficient frontier lines in the 3D space in Fig.~\ref{fig_BKE_frontier_3D}. Lines in the $ (\xi y) $-plane labeled by numbers and colors are efficient frontiers of the planning problem at the confidence levels $ \alpha $ marked on the curve. For any fixed value of $ \alpha $, any point on a line is as good as any other point in terms of the quality of the solution obtained with the model. Further preferences may be based on other considerations that were not used as a part of model inputs. 
} 
\label{fig_frontier_2D}
\end{center}
\end{figure}

We will refer to such contour lines at different probability levels as {\it efficient frontiers}, similarly to the use of this term in the Markowitz (1956) portfolio theory \cite{Markowitz}. The 3D plot of the BKE tail probability $ p(y,0| \xi) $ as a function of both controls $ y $ and $ \xi $ with efficient frontier lines is shown in Fig.~\ref{fig_BKE_frontier_3D}. The corresponding 2D view that might be easier for analysis is shown in Fig.~\ref{fig_frontier_2D}.

Efficient frontier lines shown in Figs. \ref{fig_BKE_frontier_3D} and \ref{fig_frontier_2D} can be used for visualization of {\it all} possible optimal policies at different confidence levels $ \alpha $ for a given plan with user-defined initial and target wealth levels and constraints on admissible values of control parameters. 
Note that while our original objective was formulated as an optimization problem, the solution presented here reduces to a combination of a 2D cubic spline and 
a root-finding, or equivalently a contour construction algorithm.\footnote{The Python package ContourPy (https://contourpy.readthedocs.io/en/v1.1.0/) can be useful to this end. ContourPy is incorporated within the Python package Matplotlib, which was used to produce Figs.~\ref{fig_BKE_frontier_3D} and 
\ref{fig_frontier_2D}.}  The plan optimization part thus reduces to the task of drawing a constant-probability line on a computed solution surface given by the BKE density $ p(y,t| \xi) $.

\section{Summary}
\label{sect_Summary}

This paper addresses a problem that is routinely solved (or at least supposed to be routinely solved) by millions of working individuals who participate in retirement saving plans such as the 401(K) plans in the US, or similar programs worldwide. The problem is to plan an optimal contribution schedule to their investment portfolio so that their wealth by the time of retirement will be above a certain target wealth level at a certain probability confidence level.

While his problem appears to be entirely unrelated to any sort of quantum physics, what we proposed in this paper is that in fact, {\it this problem is quantum mechanics in disguise}!
This of course does not mean quantization of the retirement wealth, or (unfortunately!) a teleportation to a wealthier quantum state, or anything of this sort that can come to the minds of non-physicists when they hear the words 'quantum mechanics' or 'Schr{\"o}dinger equation'. What we mean by saying ``retirement planning {\it is} quantum mechanics in disguise" rather refers to the fact that, with standard simplifications and assumptions of normal equity returns etc., the first problem is {\it mathematically identical} to a very famous problem in quantum mechanics, namely quantum mechanics with the Morse potential. 

It is interesting to note here that the non-linear Morse potential
(\ref{U_Morse}) arises in our approach due to the mere presence of controls $(u_0, \xi) $. Indeed, if we take $ u_t = 0 $ in Eq.(\ref{SDE_port_Pi}), we end up with the standard Geometric Brownian Motion (GBM) model as a model of the portfolio wealth process. This can be contrasted with scenarios where the 
{\it uncontrolled} dynamics are non-linear due to a non-linear potential. Models of non-linear dynamics for equity stock markets with nonlinearity generated by market flows and market frictions have been previously proposed by the author in (Halperin, 2021, 2022 \cite{IH_2021, IH_2022}.) In contrast, in this paper we start with {\it linear} uncontrolled dynamics of the GBM model, while non-linearities arise in the mathematically equivalent formulation of the problem based on a transition to the Morse variable $ y $ that blends together the state and action variables.

Our method can be compared with the conventional feedback-loop control approach of dynamic programming (DP) or reinforcement learning (RL). It is well known that low dimensional problems with continuous controls can be efficiently solved by discretization of action variables. Therefore, the mere fact of reducing our 2D control problem to a 2D grid search problem may be hardly surprising on its own. Nonetheless, for practical applications that we have in mind for this paper, the end criteria of success are the accuracy and the speed of calculation. The approach of this paper can be viewed as a way to efficiently construct a grid of policies and value functions (i.e. BKE tail probabilities, in our case) using methods that work in physics. 

While our step of constructing a grid of controls and output BKE tail probabilities is similar to the procedure employed in DP and RL upon discretization of action variable, one important difference of our scheme from DP or RL is what we do {\it next} with this grid. With DP or RL, we want to find a policy (i.e. a grid point) where the value function attains its global maximum - which can be performed via a regular grid search because of discretization. In contrast, with goal-based wealth management, our setting fits the framework of satisficing decision-making of Herbert Simon (1955, 1979) \cite{Simon_1955, Simon_1979}. This is because with goal-based wealth management, we want to find a policy that is {\it good enough} (i.e. {\it satisficing}) to drive the probability of the terminal portfolio wealth falling below some target value $ \hat{\Pi} $ to the level of $ \alpha $, but not below! This is {\it not} a maximization problem that is usually solved by DP or RL, this is a problem 
of finding {\it contour lines} on a surface of the BKE tail probability. Therefore, the nature of a solution for DP or RL is quite different from a solution to satisficing decision making: while the solution is given by a point in an action space for the former, it has a form of a degenerate contour line for the latter.     

One important concept that we propose here as a way to communicate the model results to end users is the notion of {\it efficient frontiers} - constant probability level curves on the surface of the BKE tail probability $ p(y,0| \xi ) $ viewed as a function of control parameters $ y $ and $ \xi $. Note here that most of RL methods typically provide a single point in a multi-dimensional action space as an optimal solution to a stochastic control problem. While the present paper emphasizes the importance of continuous degenerate families of solutions (efficient frontiers) for our particular 2D control problem with a satisficing  objective function, 
the author is unaware of any systematic discussion of degenerate satisficing solutions in a general RL setting.  

Back to the mathematical identity of the retirement planning problem to the quantum mechanical problem, note that we reduced our initial problem not just to the standard textbook problem of solving the Schr{\"o}dinger equation with a fixed Morse potential, but rather to a {\it quantum optimal control problem} with the Morse potential. Unlike the passive dynamics problem of the conventional quantum mechanics, in tasks of quantum control, the objective is to control the Schr{\"o}dinger potential in order to bring a quantum system into a desired terminal state at a smallest cost or in a shortest time.  In the same way, here we control the degree of non-linearity of the Morse potential and the initial particle position in order to achieve a desired quantum mechanical terminal state $ \bar{\Psi}(y,T) $. In this sense, the proposition formulated in the beginning of this section can be refined to the following proposition: \\

{\it The wealth management problem is a quantum optimal control problem in disguise}. \\
\\
where, as before, by saying {\it in disguise} we mean {\it mathematically equivalent}.
A recent work applied quantum computing to quantum mechanics of molecules with the Morse potential to compute energies of two lowest eigenstates by Apanasevich {\it et. al.} (2021) \cite{Apanasevich_2021}. Exploring potential applications of quantum computing to the problem formulated in this paper and its potential extensions can be an interesting venue for future research. 

Another interesting direction would be to explore alternative methods for the FPE and BKE equations which may or may not make connections to the Schr{\"o}dinger equation as we did in this paper. While here we solved the BKE equation and its associated Schrodinger equation using the standard eigenvalue decomposition method coupled with a non-standard way of constructing the basis, this is clearly not the only available, and likely not the most efficient way to solve these equations. In particular, it would be interesting to apply the integral transform approach of Ref.~\cite{Itkin} in our setting of optimal control. Finally, extensions to multiple dimensions can be another interesting topic for future research along the lines proposed in this paper.

\def\thesection{A}	
\setcounter{equation}{0}
\def\theequation{\thesection.\arabic{equation}}

\section*{Appendix A: Schr{\"o}dinger equation with the Morse potential}

In this appendix, we provide semi-analytical solutions for both Schr{\"o}dinger equations (\ref{SE}) and (\ref{SE_back}) which we referred to in the main text as 
FP-SE and BK-SE, respectively. As both SEs share the same Hamiltonian and only differ in the initial conditions, solutions for both equations will be found using 
a variant of the eigenvalue decomposition method, where the only difference between Eqs.(\ref{SE}) and (\ref{SE_back}) amounts to different initial conditions.
We will start our analysis with the FP-SE (\ref{SE}) and then shown where the final formulae for the BK-SE (\ref{SE_back}) differ from the final formulae for the FP-SE (\ref{SE}).

Our approach largely follows Ref.\cite{Molnar_2002} that offers a semi-analytical approach to the time-dependent SE with the Morse potential using the standard basis function decomposition method for PDEs. The novelty of this method is in the way of constructing the basis itself.
The method of \cite{Molnar_2002} gives rise to a 'natural' infinite and countable basis for the Hilbert space corresponding to the SE (\ref{SE}). 

The method of \cite{Molnar_2002} is based on a recursive algebraic procedure which is rooted in supersymmetric quantum mechanics (SUSY QM) of Witten \cite{Witten_1982}, see e.g. \cite{Cooper_SUSY_book} for a review. To ensure that the present paper is reasonably self-contained and can be read by a diverse audience that does not only include theoretical physicists, this appendix provides a short overview of this approach, together with its application to our problem.  

We start with Eqs.(\ref{Schrodinger_change}) which we repeat here for convenience: 
\bea
\label{Schrodinger_change_A}
f(x,t) &\hspace{-0.2cm} = \hspace{-0.2cm} &  e^{ - \frac{V(x)}{\hbar} } \Psi(x,t) \nonumber \\
p(x,t) &\hspace{-0.2cm} = \hspace{-0.2cm} & e^{  \frac{V(x)}{\hbar} } \bar{\Psi}(x,t)
\eea
The Schr{\"o}dinger equations (\ref{SE}) and (\ref{SE_back}) read
\beq
\label{SE_A_0}
- \hbar \frac{\partial \Psi(x,t)}{\partial t} = \mathcal{H}_{\pm} \Psi(x,t)
\eeq
\beq
\label{SE_back_A}
- \hbar \frac{\partial \bar{\Psi}(x,\tau)}{\partial \tau} = \mathcal{H}_{\pm} \bar{\Psi}(x,\tau)
\eeq
where we use the Hamiltonians $ \mathcal{H}_{+} $ or $ \mathcal{H}_{-} $ when, respectively, $ g > 0 $ or $ g < 0 $. The latter are defined as follows:
\beq
\label{H_plus}
\mathcal{H}_{+}:= - \frac{\hbar^2}{2} \frac{\partial^2}{\partial x^2} + D_0 \left(e^{-2x} - 2 e^{-x} + 1 \right) + U_0, \; \; \; \left(g > 0 \right) 
\eeq
\beq
\label{H_minus}
\mathcal{H}_{-} := - \frac{\hbar^2}{2} \frac{\partial^2}{\partial x^2} + D_0 \left(e^{-2x} + 2 e^{-x} + 1 \right) + U_0, \; \; \; \left( g < 0 \right)
\eeq
Focusing on the first Hamiltonian $ \mathcal{H}_{+} $, it is the same as the Hamiltonian $ \hat{H}_0 $ in \cite{Molnar_2002}, with $ \alpha = m = 1 $ and $ \hat{P} = - i \frac{\partial}{\partial x} $. To keep consistency of our notations with \cite{Molnar_2002}, we introduce 
parameter $ s $ as follows:
\beq
\label{omega_0_s}
s := \frac{\sqrt{2 D_0}}{\hbar} - \frac{1}{2} = \frac{1}{2} \left(g-1\right) = \frac{\eta}{\hbar}  
\eeq

The Hamiltonian $ \mathcal{H}_{+} $ can be expressed in terms of a new dimensionless Hamiltonian $ H_{+} $:
\beq
\label{H_0}
\mathcal{H}_{+}  = 
\frac{\hbar^2}{2} H_{+},
\; \; \; H_{+} := \hat{P}^2 + \left(s + \frac{1}{2} \right)^2 \left(e^{-2x} - 2 e^{-x} + 1 \right) + E_0, \; \; \; 
E_0 := \frac{2 U_0}{\hbar^2}
\eeq  
Similarly, for the second Hamiltonian (\ref{H_minus}), we obtain:  
\beq
\label{tilde_H_0}
\mathcal{H}_{-} = 
\frac{\hbar^2}{2} H_{-},
\; \; \; H_{-} := \hat{P}^2 + \left(s + \frac{1}{2} \right)^2 \left(e^{-2x} + 2 e^{-x} + 1 \right) + E_0  
\eeq  
Introducing the dimensionless time parameter
\beq
\label{tau}
\tau :=  \frac{\hbar}{2}  t
\eeq
(not to be confused with $ \tau := T - t $ that appears in the BK-SE (\ref{SE_back_A})), the SE (\ref{SE_A_0}) is written in dimensionless variables $ (x, \tau) $ as follows:
\beq
\label{SE_A}
-  \frac{\partial \Psi(x, \tau)}{\partial \tau} = H_{\pm} \Psi(x, \tau)
\eeq
If $ \Psi(x, \tau) $ is a solution of (\ref{SE_A}), then the solution of the original SE (\ref{SE_A_0}) is $ \Psi \left(x, \hbar t/2 \right) $. Similarly, the BK-SE in dimensionless units $ (x,\tau) $ reads
\beq
\label{SE_A_plus}
-  \frac{\partial \bar{\Psi}(x, \tau)}{\partial \tau} = H_{\pm} \bar{\Psi}(x, \tau)
\eeq
Next we define the so-called generalized SUSY ladder operators
\bea
\label{ladder_ops}
\mathcal{Q}(q) &\hspace{-0.2cm} := \hspace{-0.2cm} & q - \left(s + \frac{1}{2} \right) e^{-x} + \frac{\partial}{\partial x} \nonumber \\
\mathcal{Q}^{+}(q) &\hspace{-0.2cm} := \hspace{-0.2cm} & q - \left(s + \frac{1}{2} \right) e^{-x} - \frac{\partial}{\partial x}
\eea  
where $ q $ is a parameter to be specified below.
 These operators satisfy the commutation relation
\beq
\label{commutator}
\left[\mathcal{Q}(q), \mathcal{Q}^{+}(q') \right] = - \mathcal{Q}(q) - \mathcal{Q}^{+}(q') + q + q'
\eeq
If we take $ q = s $ in (\ref{ladder_ops}), the Hamiltonian $ H_{+} $ factorizes in terms of ladder operators:
\beq
\label{H_0_fact}
H_{+} = \mathcal{Q}^{+}(s)\mathcal{Q}(s) + \varepsilon_0, \; \; \; \varepsilon_0 :=  \frac{2 U_0}{\hbar^2} + s + \frac{1}{4}
\eeq
The second Hamiltonian $ H_{-} $ can be factorized in terms of the same operators 
as follows: 
\beq
\label{H_0_fact_2}
H_{-} = \left(2s - \mathcal{Q}(s) \right)  \left(2s - \mathcal{Q}^{+}(s) \right) + \varepsilon_0 
\eeq
Let us again focus for now on the first Hamiltonian $ H_{+} $ which has a minimum. A lowest energy state (a ground state) $ \left| \phi_0 \rangle \right. $ 
of $ H_{+} $ with energy equal to $ \varepsilon_0 $ could be obtained, if it exists, from the requirement that it should be annihilated by the ladder operator $ \mathcal{Q}(s) $:
\beq
\label{ground_state}
\mathcal{Q}(s)  \left| \phi_0 \rangle \right. = 0 
\eeq
Using the explicit expression for $ \mathcal{Q}(s) $ from Eqs.(\ref{ladder_ops}) and solving the resulting ODE, we obtain the lowest energy eigenstate  
of the Morse potential in the coordinate representation, which we denote $ \phi_0(x) $ and call the ground state wave function (WF):
\beq
\label{ground_state_x}
\phi_0(x) = c_0 e^{- s x - \left(s + \frac{1}{2} \right)e^{-x} }
\eeq    
where $ c_0 $ is a normalization constant that should be found from the normalization condition $ \int \phi_0^2(x) dx = 1 $, provided this integral 
converges.\footnote{
It is only in the latter case when $ \phi_0(x) $ is normalizable that it really exists as a true ground state.
In certain other cases encountered in SUSY QM, even though non-normalizable analytical solutions to 
Eq.(\ref{ground_state}) may be available, their corresponding normalization integrals may diverge, indicating that no stationary ground states exist for such problems.}
As is seen from (\ref{ground_state_x}), the ground state WF (\ref{ground_state_x})  is normalizable only if $ s > 0 $. With this choice, the ground state WF (\ref{ground_state_x}) decays exponentially for $ x \rightarrow + \infty $ and super-exponentially for $ x \rightarrow - \infty $.

Factorization properties of certain Hamiltonians similar to Eq.(\ref{H_0_fact}) play a critical role in supersymmetric quantum mechanics (SUSY QM) of Witten 
\cite{Witten_1982, Cooper_SUSY_book}, where they are used to build higher excited states algebraically by a sequential application of ladder operators (\ref{ladder_ops}) to lower-energy states. However, such direct methods of building higher-energy states are not particularly useful for the problem of solving the time-dependent SE (\ref{SE_A}) with the Morse potential. This is because the Morse potential has only  $[s]+1 $ bound states when $ s > 0 $ (where $ [s]$ is the largest integer that is still smaller than $ s $) with normalizable discrete-energy wave functions (WFs) \cite{Landau}. Furthermore, even when $ s > 0 $, a finite set of normalizable WFs does {\bf not} form a complete orthogonal basis in the infinite-dimensional Hilbert space corresponding to the QM problem (\ref{SE_A}) \cite{Molnar_2002}.\footnote{In our numerical experiments, typical values of inputs from the planner and market produced values of $ s $ in the range $[-1/2,1/2] $,
which produces at most one bound state when $ 0 < s < 1/2 $.}   

To conclude so far, the standard SUSY approach of constructing a convenient basis is not applicable for the present problem, either because no stable ground state exists (as in happens when $ s < 0 $) or because higher states are not normalizable (which happens when $ s > 0 $).
   
Instead of using a ground state WF (\ref{ground_state_x}) (which may not even exist if  $ s < 0 $) as a starting state to build a basis, 
we consider the following definition of a {\it new} starting state $ | \phi_0 \rangle $ \cite{Molnar_2002}:  
%
\beq
\label{ground_state_q}
\mathcal{Q}(q)  \left| \phi_0 \rangle \right. = 0 
\eeq
Here $ q $ is a parameter that will be specified below. This produces the following normalized state in the coordinate representation: 
\beq
\label{ground_state_x_q}
\phi_0(x) = \frac{(2s+1)^q}{\sqrt{\Gamma(2q)}}  e^{- q x - \left(s + \frac{1}{2} \right)e^{-x} }
\eeq    
This shows that the new state $ | \phi_0 \rangle $ is normalizable as long as $ q > 0 $ and $ s > - 1/2 $. 



Using the initial state (\ref{ground_state_x_q}), an infinite set of orthogonal squared-normalized states is now constructed as follows \cite{Molnar_2002}:
\beq
\label{basis_Molnar_phi}
\left| \phi_n \right\rangle  = \left\{  \prod_{k=1}^{n} C_{k}^{-1} \mathcal{Q}^{+}(q + k -1) \right\} \left. \left. \right| \phi_0 \right \rangle, \; \; \; n = 1, \ldots, \infty
\eeq
where coefficients $ C_{k} = \sqrt{k(k+ 2q - 1)} $ enforce the correction normalization. States (\ref{basis_Molnar_phi}) are referred to 
as {\it quasi-number states}. Eqs.(\ref{basis_Molnar_phi}) and (\ref{commutator}) imply that the ladder operators act on the quasi-number states as raising and lowering operators:
\bea
\label{raising_lowering}
\mathcal{Q}^{+} (q + n) \left| \phi_n \right\rangle &\hspace{-0.2cm} = \hspace{-0.2cm} & C_{n+1} \left| \phi_{n+1} \right \rangle \nonumber \\
 \mathcal{Q} (q + n) \left| \phi_n \right\rangle &\hspace{-0.2cm} = \hspace{-0.2cm} & C_{n} \left| \phi_{n-1} \right \rangle, \; \; \; n = 0,1, \ldots
\eea
where $ C_0 = 0 $ and $ \left| \phi_{-1} \right \rangle = 0 $. On the other hand, these relations imply that the quasi-number state $  \left| \phi_{n} \right \rangle  $ is an eigenstate of the operator 
$ \mathcal{Q}^{+} (q + n-1) \mathcal{Q} (q + n) $:
\beq
\label{quasi_number_as_eigen}
\mathcal{Q}^{+} (q + n-1) \mathcal{Q} (q + n) \left| \phi_{n} \right \rangle = C_n^2 \left| \phi_{n} \right \rangle
\eeq
As shown in \cite{Molnar_2002}, in the coordinate representation 
Eqs.(\ref{basis_Molnar_phi}) produce the following expressions:
\beq
\label{basis_Molnar}
\phi_n(x) = \sqrt{\frac{n!}{\Gamma(2q + n)}} y^{q}(x) \exp\left( - \frac{y(x)}{2} \right) L_{n}^{(2q -1)} \left(y(x) \right), \; \; \; n = 0,1, \ldots
\eeq  
where $ y(x) := (2 s + 1) e^{-x} $ is the so-called Morse variable (so that $ 0 \leq y(x) \leq \infty $), and $ L_{n}^{(2q -1)}(y) $ are generalized Laguerre polynomials 
\cite{AS}. 
Therefore, basis functions $ \phi_n(x) $ depend on $ x $ only through the Morse variable $ y(x) $. In contrast to energy eigenfunctions, the quasi-number states  
(\ref{basis_Molnar}) form a complete orthonormal set of square integrable functions, and thus can be used as a basis in the infinite dimensional Hilbert space 
corresponding to the SE (\ref{SE_A})  \cite{Molnar_2002}. 

\begin{figure}[ht]
\begin{center}
\hspace*{-1cm} \includegraphics[
width=175mm,
height=65mm]{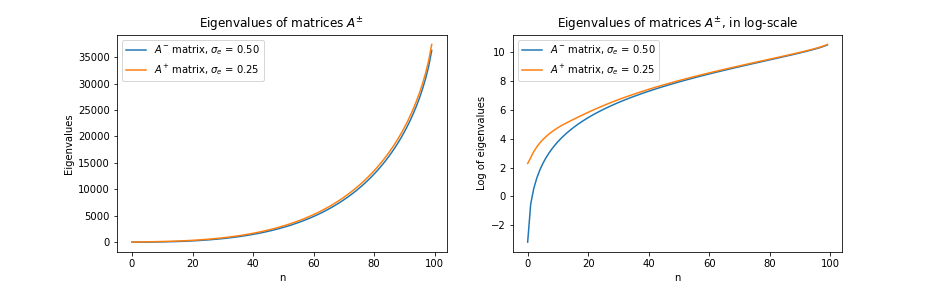}
\caption{Eigenvalues of matrices $ {\bf A}^{\pm} $ truncated at size $ N = 100 $. On the left: Eigenvalues grow quadratically with $ n $ for large values of $ n $. On the right: eigenvalues are displayed in log-scale to show differences in spectra for $ {\bf A}^{+} $ and $ {\bf A}^{-} $.
} 
\label{fig_eigenvals}
\end{center}
\end{figure}

Next we use the orthogonal basis (\ref{basis_Molnar}) to compute matrix elements of the Hamiltonians with this basis. Let  $ {\bf A}^{\pm} $ be matrices with matrix elements $ A_{mn}^{\pm} = \langle \phi_m | H_{\pm} | \phi_n \rangle $. They can be found  
using Eqs.(\ref{commutator}), (\ref{H_0_fact}) and (\ref{H_0_fact_2}) as follows: 
\bea
\label{H_mn}
A_{mn}^{+} := \langle \phi_m | H_{+} | \phi_n \rangle 
&\hspace{-0.2cm} = \hspace{-0.2cm} & \left( C_n^2 + \left(s-q -n \right)^2 + \varepsilon_0 \right) \delta_{m,n} \nonumber \\ 
&\hspace{-0.2cm} + &\hspace{-0.2cm}
\left(s -q - n \right) C_m \delta_{m,n+1} +  \left(s-q - m \right) C_n \delta_{m,n-1}  \nonumber \\
A_{mn}^{-} := \langle \phi_m | H_{-} | \phi_n \rangle &\hspace{-0.2cm} = \hspace{-0.2cm} &
\left( C_n^2 + \left(s - q - n \right)^2 + 2(n+q)(2s+1) + \varepsilon_0 \right) \delta_{m,n}  \\
&\hspace{-0.2cm} - \hspace{-0.2cm} &  
\left(s + q + n +1 \right) C_m \delta_{m,n+1} - \left( s + q + m +1 \right) C_n \delta_{m+1,n}  \nonumber
\eea 
Eqs.(\ref{H_mn}) show that $ {\bf A}^{\pm} $ are infinite-dimensional symmetric tridiagonal matrices. In practical applications, the basis is truncated at some value $ N $, so that 
we end up with an $ N$-dimensional truncated version of the original matrices $ {\bf A}^{\pm} $.\footnote{In \cite{Molnar}, the tridiagonal form of the generator 
$ {\bf A}^{+} $ was used to take the value $ q = s - [s] $ where $ [s] $ in the maximal integer that is still smaller than $ s $, to decouple bounded states from unbounded states from the continuous spectrum. This prescription assumes that $ s . 0 $, and therefore cannot be applied in our setting when $ s < 0 $.}  Such tridiagonal matrices can diagonalized by 
orthogonal transformations:
\beq
\label{A_eigenvalues}
{\bf A}^{\pm} = {\bf U}_{\pm}\boldcalD_{\pm} {\bf U}_{\pm}^T
\eeq
where $ \boldcalD_{\pm} $ are diagonal matrices of eigenvalues, and $ {\bf U}_{\pm} $ are orthogonal matrices of eigenvectors that satisfiy the relations $ {\bf U}_{\pm} {\bf U}_{\pm}^T = 
 {\bf U}_{\pm}^T {\bf U}_{\pm} = \mathbb{I} $. Eigenvalues of matrices $ {\bf A}^{\pm} $ exhibit a quadratic dependence 
 on $ n $ for large $ n $, see Fig.~\ref{fig_eigenvals}. One can notice that such quadratic dependence is similar to one observed for the Schr{\"o}dinger equation (\ref{SE_A}) for {\it bounded} states \cite{Landau}. Recall that bounded states exist for the Morse potential only when $ s > 0 $, and there are $ [s]+1 $ of them. While in our case
 $ s < 0 $ and Hamiltonians $ H_{\pm} $ do not produce bound states, the eigenvalue spectra for both $ H_{+} $ and $ H_{-} $ in the orthogonal basis 
 (\ref{basis_Molnar}) have 
 a quadratic dependence on $ n $ for arbitrary large values of $ n $.
 
 The method of \cite{Molnar_2002} thus elegantly uses the SUSY structure of the SE (\ref{SE_A}) to come up with a `natural' countable basis for this problem.
To proceed, we follow standard textbook steps, and 
look a solution in the form of an eigenstate decomposition with the basis (\ref{basis_Molnar}) and time-dependent coefficients $ w_N^{(F)}(t) $:
\beq
\label{eigen_val_exp}
\Psi(x,\tau) = \sum_{n=0}^{\infty} w_n^{(F)}(\tau) \phi_n(x)
\eeq
Here we denoted weights as $ w_n^{(F)}(t) $ to differentiate them from weights $ w_N^{(B)}(t) $ that will be introduced below for the BK-SE (\ref{SE_back}).
We first use this relation to find the initial weights $ w_n^{(F)}(0) $ at time 0. They can be found from the initial condition 
\beq
\label{init_cond_Psi_A}
 \Psi(x,0) = \delta(x - x_0) e^{\frac{V(x_0 + \hat{z})}{\hbar}}
 \eeq  
which follows from Eq.{(\ref{init_cond_Psi}), and the completeness relation for basis functions
\beq
\label{orthogo}
\sum_{n=0}^{\infty} \phi_n (x) \phi_n(x') = \delta( x - x')
\eeq
This produces the initial weights
\beq
\label{init_weight}
w_n^{(F)}(0) = \int_{-\infty}^{\infty} dx \Psi(x,0) \phi_n(x) = \int_{-\infty}^{\infty} dx e^{ V(x)/\hbar} \delta(x - x_0) \phi_n(x) 
= e^{V(x_0)/\hbar } \phi_n(x_0)
\eeq
Note that $ w_n^{(F)}(0) $ scale with $ n $ for $ n \rightarrow \infty $ as $ n^{-1/4} \sin \left[ 2 \sqrt{n y_0} - \pi
\left( q - \frac{3}{4} \right) \right] $,   
as can be seen using the asymptotic expression for Laguerre polynomials for large values of $ n $:
\beq
\label{Laguerre_asymptotics}
L_n^{(\alpha)}(x) = \frac{n^{\alpha/2 - 1/4}}{\sqrt{\pi}} x^{- \alpha/2 - 1/4} e^{x/2} \sin \left( 2 \sqrt{nx} - \frac{\pi}{2} \left(\alpha - 1/2 \right) \right) + 
O \left( n^{\alpha/2 - 3/4} \right)
\eeq 
A slow oscillating decay of initial weights with the basis index $ n $ is illustrated in Fig.~\ref{fig_initial_weights} on the left. Note that strong oscillations with $ n $ are only observed for initial weights. Finite-$T$ weights for long enough times $ T $ show a much smoother behavior with a fast decay of weights 
with $ n $, as shown on the right of Fig.~\ref{fig_initial_weights}. Fig.~\ref{fig_weights_FPE_all_T} shows the behavior of finite-$T$ weights $ w_n^{(F)}(T) $ 
as functions of  $n $
for different planning horizons $ T $, and shows that their decay with $ n $ for sufficiently large values of $ T $ such as $ T = 15 $ or $ T = 20 $ is exponential in $ n $.
 
  \begin{figure}[ht]
\begin{center}
\hspace*{-1cm} \includegraphics[
width=175mm,
height=65mm]{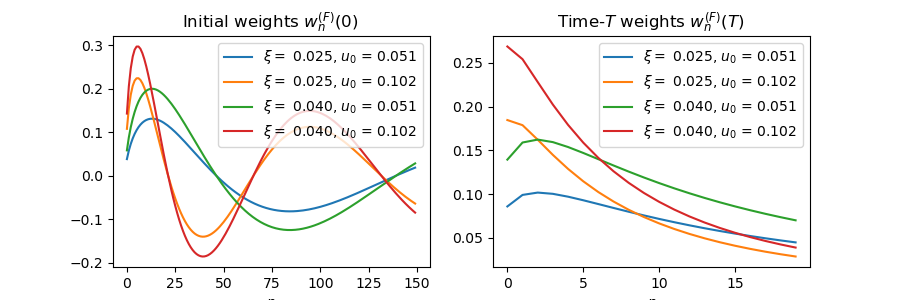}
\caption{On the left: Initial weights $ w_n^{(F)}(0) $ in Eq.(\ref{init_weight}) as a function of $ n $ for different values of $ u_0 $. On the right: the finite-$ T $ coefficients $ w_n{(F)}(T) $ for $ T=20 $. 
} 
\label{fig_initial_weights}
\end{center}
\end{figure}

\begin{figure}[ht]
\begin{center}
\hspace*{-1cm} \includegraphics[
width=185mm,
height=65mm]{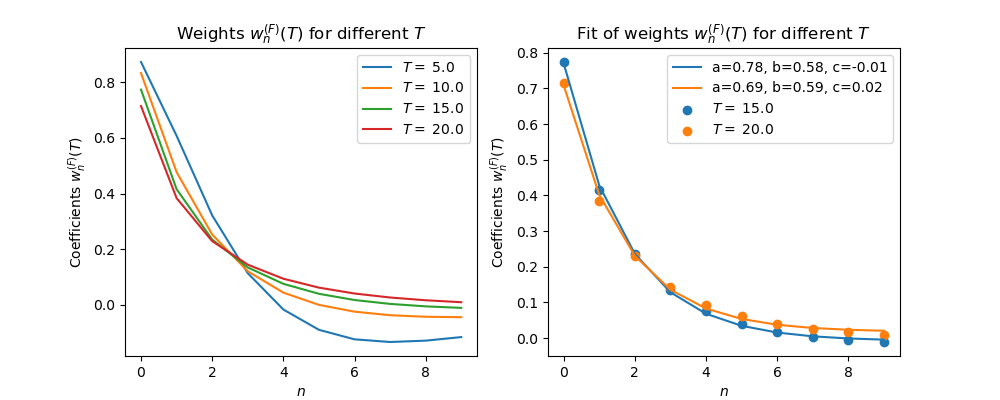}
\caption{On the left: finite-$T $ weights $ w_n(T) $ for the FP-SE (\ref{SE}) as functions of $ n $ for different planning horizons $ T $. Unlike initial weights $ w_n(0)$ that decay as $ n^{-1/4} $ with additional oscillations (see Fig.~\ref{fig_initial_weights}) on the left), weights for 
long times exhibit a fast and smooth decay with $ n $. On the right: fit of weights $ w_n(T) $ for $ T=15 $ and $ T=20 $ by the function $ \hat{w}_n = a e^{- b n} + c $ with parameters $ a, b, c$.  
} 
\label{fig_weights_FPE_all_T}
\end{center}
\end{figure}

Still, while finite-$ T$ coefficients $ w_n^{F)} $ show a smooth regular behavior with a fast decay with  $ n $, the slow decay of initial weights $ w_n^{F)}(0) $ 
might be concerning and pointing to a possible loss of accuracy for the FPE equation due to truncation of the original infinite-dimensional basis.
It is therefore of interest to analytically explore implications of truncating the basis to a finite value of $ N $ such as $ N = 100 $ or $ N = 150 $ for numerical implementation. 

To begin, we start with an observation that a finite-$N$ approximation to the FPE density $ f(x,0) $ is no longer exactly equal to $ \delta(x - x_0) $. Instead, it can be computed in closed form in terms of the Morse 
variable $ y = y(x) = (2s+1)e^{-x} $ using Eq.(10.12.9) in 
\cite{BE}:
\beq
\label{sum_L_n_L_n}
\sum_{n=0}^{N-1} \frac{n! L_n^{(2q-1)}(x) L_n^{(2q-1)}(y)} {\Gamma(n+2q)}  = \frac{N!}{\Gamma(N+2q-1)} \frac{1}{x-y} \left[ L_{N-1}^{(2q-1)}(x) 
L_N^{(2q-1)}(y) - L_N^{(2q-1)}(x) L_{N-1}^{(2q-1)}(y) \right]
\eeq
In the limit  $ N \rightarrow \infty $, this expression can be approximated as follows using the asymptotic form of $ L_n^{(2q-1)}(x) $:
\bea
\label{sum_LN_2}
 \sum_{n=0}^{N-1} \frac{n! L_n^{(2q-1)}(x) L_n^{(2q-1)}(y)}{\Gamma(n+2q)}  
&\hspace{-0.2cm} = \hspace{-0.2cm}  & \frac{(x y)^{-q+1/4} 
e^{(x+y)/2}}{2 \pi(x-y)}  \left[ 
  \left( \sqrt{x} + \sqrt{y} \right) \sin\left(2 \sqrt{N}(\sqrt{x} - \sqrt{y} ) \right) \right. \nonumber  \\
&\hspace{-0.2cm} - \hspace{-0.2cm} &
  \left. \left( \sqrt{x} - \sqrt{y} \right) \sin\left(2 \sqrt{N}(\sqrt{x} + \sqrt{y}) - 
 2 \pi (q - 3/4)  \right)  
 \right]  
\eea
The last relation can be used to explore the behavior of a finite-$N$ approximation to the FPE density $ f(x,0) $ at time zero, which we will denote $ f_N(x,0) $.
We obtain
\bea
\label{FPE_finite_N}
f_N(x,0) 
&\hspace{-0.2cm} = \hspace{-0.2cm} & \frac{1}{2 \pi} y^{-1/2} \left(\frac{y}{y_0} \right)^{-\eta/ \hbar - 1/4} e^{ - \frac{y-y_0}{2}} \left[ 
\frac{ \sin \left[ 2 \sqrt{N} \left( \sqrt{y} - \sqrt{y_0} \right) \right]}{ \sqrt{y} - \sqrt{y_0} } \right.  \\
&\hspace{-0.2cm} + \hspace{-0.2cm} & \left. \sin (2 \pi q) \frac{ \sin \left[ 2 \sqrt{N} \left( \sqrt{y} + \sqrt{y_0} \right) \right]}{ \sqrt{y} + \sqrt{y_0}} 
+ \cos (2 \pi q) \frac{ \cos \left[ 2 \sqrt{N} \left( \sqrt{y} + \sqrt{y_0} \right) \right]}{ \sqrt{y} + \sqrt{y_0}} \right] \nonumber
\eea
We now can use this expression to guide our choice for parameter $ q > 0 $ which was considered a free parameter up to this point. In the strict limit
$ N \rightarrow \infty $, the first two terms in this expression converge, respectively, to $ \delta(\sqrt{y} - \sqrt{y_0}) $ and $ \delta(\sqrt{y} + \sqrt{y_0}) $, and furthermore the second delta-function $ \delta(\sqrt{y} + \sqrt{y_0}) $ vanishes for all values of $ y_0, y \neq 0 $. On the other hand, 
the third term does not have a well-defined limit $ N \rightarrow \infty $ for generic values of $ q $. Therefore, we pick a value of $ q $ that satisfies the relation $ 2 \pi q = m \pi/2 $ with $ m = 1, 3, 5, \ldots $, or equivalently 
\beq
\label{q_choices}
q = \frac{m}{4}, \; \; \; \; m = 1,3,5, \ldots
\eeq 
Furthermore, analysis of the probability current in this model performed in Appendix C shows that parameter $ q $ should be larger than $ \eta/ \hbar $.
As in our practical setting the ratio $ \eta/\hbar $ varies approximately between -0.2 and 0.2, our final choice for $ q $ corresponds to (\ref{q_choices}) with $ m = 3 $:
\beq
\label{q_final}
q = \frac{5}{4}
\eeq
Assuming this choice in (\ref{FPE_finite_N}) and transforming the delta-function $ \delta(\sqrt{y} - \sqrt{y_0} ) $ back from the Morse variable $ y = y(x) $ to the original $x$-space, we recover the correct initial condition in the limit  $ N \rightarrow \infty $:
 \beq
 \label{init_f_from_finite_N_limit}
 \lim_{N \rightarrow \infty} f_N(x,0)  = \delta(x - x_0) 
 \eeq
 On the other hand, we see that for a finite $N$ approximation to the initial density $ f(x,0) $ produces strong oscillations around the singular point $ y = y_0 $, see Fig.~\ref{fig_finite_N_FPE_t0}. 
 
\begin{figure}[ht]
\begin{center}
\hspace*{-1cm} \includegraphics[
width=80mm,
height=65mm]{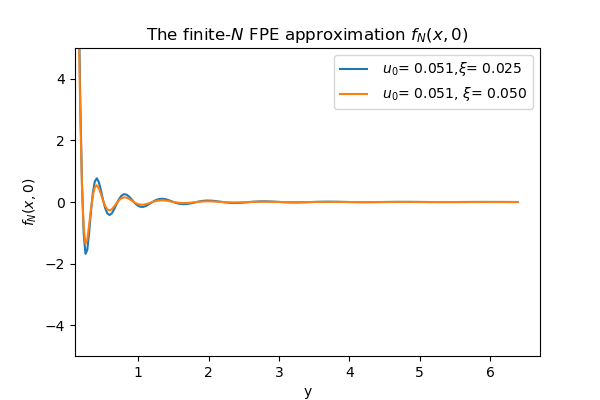}
\caption{The finite-$N$ approximation for the FPE density $ f_N(y,0) $ at $ t = 0 $ with $ N = 150 $, computed using Eq.(\ref{FPE_finite_N}) with $q = 5/4 $. Oscillations around the singular point $ y = y_0 $ are due to the truncation of the original infinite-dimensional basis.The correct initial state is proportional to $ \delta(y-y_0) $.
} 
\label{fig_finite_N_FPE_t0}
\end{center}
\end{figure}

Now let us consider the BK-SE (\ref{SE_back}). For this equation, the only difference from the FP-SE (\ref{SE}) is in a different initial condition
\beq
\label{init_cond_bar_Psi}
 \bar{\Psi}(x,0) =  e^{-V(x)/\hbar} \theta(\hat{x} - x)
 \eeq  
which follows from (\ref{init_cond_Psi}). 
This translates into different initial weights:
\beq
\label{init_weights_back}
w_N^{(B)} (0) = \int_{-\infty}^{\infty} dx \phi_n(x) e^{- V(x)/\hbar} \theta(\hat{x} - x)
\eeq
which gives 
\beq
\label{init_weight_back_2}
w_N^{(B)}(0) = 
\sqrt{ \frac{n!}{\Gamma(n + 2q)} } \int_{\hat{y}}^{\infty} dy z^{\kappa -1} e^{-y} L_n^{(2q-1)}(y)  
\eeq
where parameters $ \hat{y} $  and $ \kappa $ are defined as follows
\beq
\label{y_hat}
 \hat{y} =  \frac{2 u_0 e^{\xi T} }{\hbar V_{\text{target}}}, \; \; \;   \kappa := q - \eta/ \hbar
 \eeq
Here $ \hat{y} $ is  the target value in terms of the Morse variable, and $ \kappa = q - \eta/ \hbar $ is an auxiliary parameter introduced to simplify formulae to follow.
 Note at this point that the integral in (\ref{init_weight_back_2}) can also be interpreted as the integral $ \int_{\hat{y}}^{\infty} w_n^{(F)}(y) dy $ - which is of course as expected based on the Feynman-Kac representation (\ref{tail_prob_BKE}) when combined with an eigenvalue decomposition method. This indicates that the method based on the BK-SE (\ref{SE_back_A}) should be preferable to the method based on the FP-SE (\ref{SE_A_0}) only if one or more of the following applies: (1) pre-asymptotic behavior for large but finite $ N $ is better controllable within the BK-SE than the FP-SE, and (2) weights (\ref{init_weight_back_2}) can be efficiently computed, preferably in an analytical form.


Starting with the second item in this list, the weights $ w_N^{(B)}(0) $ can be computed by writing the integral in (\ref{init_weight_back_2})
as the difference of two integrals on intervals $ [0, \infty]$
 and $ [0,\hat{y}] $, and then using relations (7.414.11) and (7.415) from \cite{GR}:
 \bea
 \label{integral_y_hat}
 \int_{\hat{y}}^{\infty} e^{-x} x^{\kappa -1} L_{n}^{(\alpha)}(x) dx 
&\hspace{-0.2cm} = \hspace{-0.2cm}  & \frac{\Gamma(\kappa) \Gamma(\alpha + n + 1 - \kappa)}{n! \Gamma(\alpha + 1 - \kappa)} 
 \\
 &\hspace{-0.2cm} - \hspace{-0.2cm}  &  \frac{ \Gamma(\kappa) \Gamma(\alpha + n + 1)}{n!} \hat{y}^{\kappa} \hypgeo{2}{2}(\alpha+n+1, \kappa; \alpha+1, \kappa+1; - \hat{y}) \nonumber
 \eea  
 Here $ \hypgeo{2}{2}(a_1,a_2; b_1, b_2; z) $ is the generalized hypergeometric 
function $ \hypgeo{p}{q}(a_1,\ldots,a_p; b_1, \ldots, b_q; z) $ with $ p = q = 2 $. The latter is defined, for generic values of $ p, q $, by the following infinite series \cite{AS, BE}:
 \beq
 \label{hypgeo_p_q}
 \hypgeo{p}{q}(a_1,\ldots, a_p; b_1, \ldots, b_q; z)  = \sum_{n=0}^{\infty} \frac{ \left(a_1 \right)_n \ldots \left(a_p \right)_n}{ \left( b_1 \right)_n \ldots \left( b_q \right)_n}
 \frac{z^n}{n!}
 \eeq 
where $ \left(x \right)_n :=  \Gamma(x + 1)/\Gamma(x + 1 - n ) $ is the Pochhammer symbol. Recall that for symmetric cases $ p = q $, including our case $ p = q = 2 $, the generalized hypergeometric function $ \hypgeo{p}{q}({\bf a}; {\bf b}; z) $  is an entire function in the complex plane \cite{AS}, so that the series (\ref{hypgeo_p_q}) converges 
everywhere for this case, see Fig.~\ref{fig_2F2}. This enables a simple and numerically efficient implementation of $ \hypgeo{2}{2}(a_1,a_2; b_1, b_2; z) $, available 
e.g. in the Python package mpmath (https://mpmath.org) and Mathematica (https://www.wolfram.com). 
Weights $ w_n^{(B)}(0) $ are then computed according to Eqs.(\ref{init_weight_back_2}) and (\ref{integral_y_hat}).\footnote{
A straightforward calculation of weights (\ref{init_weight_back_2}) using the definition of the Laguerre polynomials produces a series involving the incomplete gamma function $ \Gamma(s,x) $:
\[
 w_N^{(B)}(0) = \sqrt{\frac{\Gamma(n + 2q)}{n!} } \sum_{m=0}^{n} (-1)^{m} \frac{\Gamma(\kappa + m, \hat{y})}{
 m! (n-m)!}
 \]
 This expression is however not convenient for numerical implementation as it involves a highly oscillating sign-alternating sum.
}

\begin{figure}[ht]
\hspace*{-2cm} \begin{center}
\includegraphics[
width=95mm,
height=85mm]{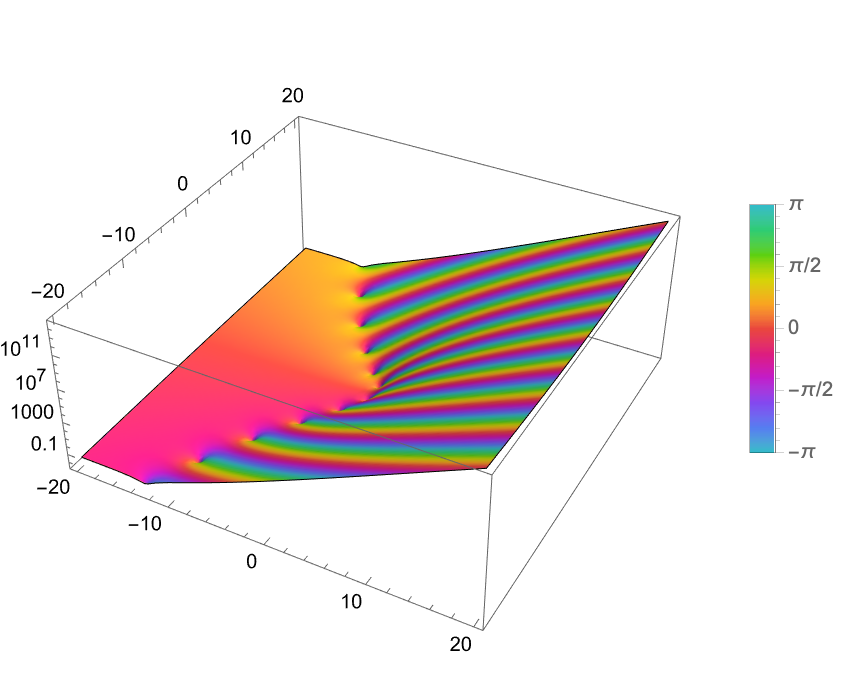}
\caption{The generalized hypergeometric function $ _2F_2(\alpha + n + 1, \kappa; \alpha+1,\kappa + 1;z) $
 in the complex $ z $-plane for 
$ \alpha = 3/2, \, \kappa = 0.7  $ and $ n = 10 $. For $ | z | \rightarrow \infty $, $ _2F_2(\alpha + n + 1, \kappa; \alpha+1,\kappa + 1;z) \sim z^{n-1}e^{z} $.  
} 
\label{fig_2F2}
\end{center}
\end{figure}

Initial weights $ w_n^{(B)}(0) $ computed according to (\ref{init_weight_back_2}) and (\ref{integral_y_hat}) are shown together with their time-$T$ values in  Fig.~\ref{fig_init_and_T_weights_BK} for $ T = 20 $. 
For comparison, in Fig.~\ref{fig_weights_BKE_all_T} we show coefficients $ w_n^{(B)}(T) $ for different planning horizons $ T $. Again, we observe an exponential decay of weights with $ n $ for sufficiently large values of $ T $, which is similar to the behavior observed for forward weights $ w_n^{(F)}(T) $ in Fig.~\ref{fig_weights_FPE_all_T}.


\begin{figure}[ht]
\begin{center}
\hspace*{-1cm} \includegraphics[
width=185mm,
height=60mm]{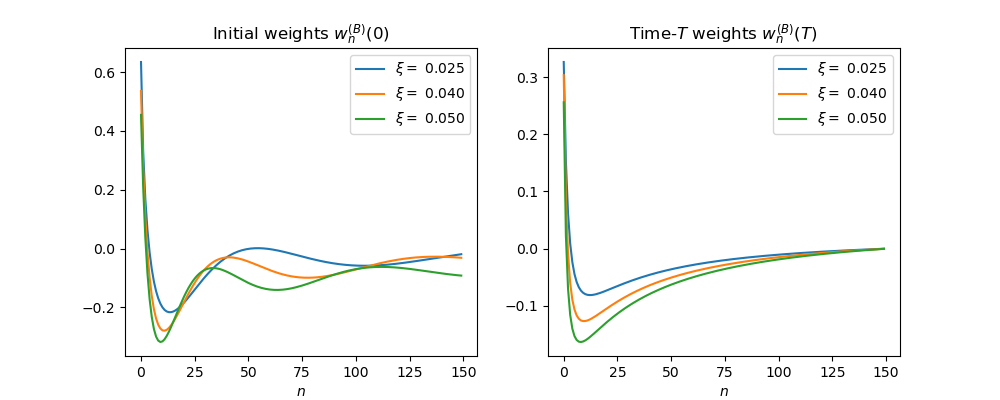}
\caption{On the left: initial weights $ w_n^{(B)}(0) $ for the BK-SE (\ref{SE_back}) as functions of $ n $ for different values of $ \xi $. On the right: time-$T $ weights 
$ w_n^{(B)}(T) $  for $ T = 20 $.  
} 
\label{fig_init_and_T_weights_BK}
\end{center}
\end{figure}

\begin{figure}[ht]
\begin{center}
\hspace*{-1cm} \includegraphics[
width=185mm,
height=60mm]{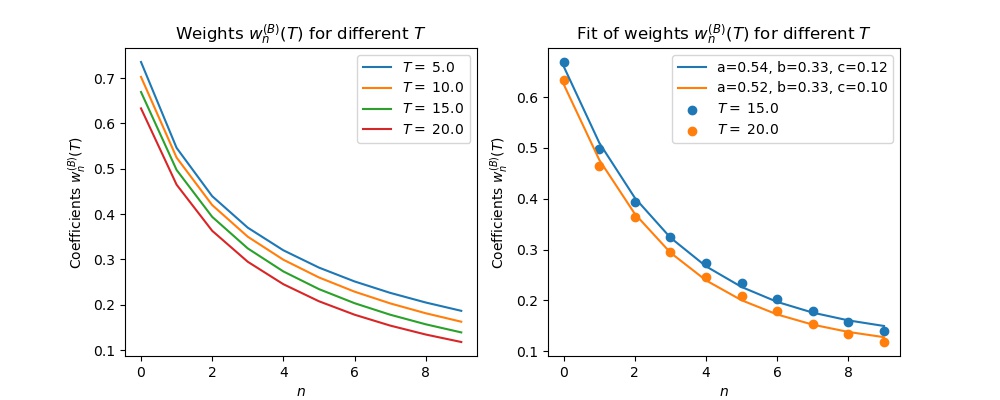}
\caption{On the left: time-$T$ weights $ w_n^{(B)}(T) $ for the BK-SE (\ref{SE_back}) as functions of $ n $ for different values of $ T $. On the right: fit of time-$T $ weights 
$ w_n^{(B)}(T) $  for $ T = 15 $ and $ T = 20 $ by function $ ae^{-bx} + c $.  
} 
\label{fig_weights_BKE_all_T}
\end{center}
\end{figure}

To explore the large-$N$ behavior of $ \bar{\Psi}_N(x,0) $, we again use (\ref{sum_LN_2}) where we now set $ q = 5/4 $ according to (\ref{q_final}). Substituting this relation  into the formula $ \bar{\Psi}(x,0) = \sum_{n} w_n^{(B)}(0) \phi_n(x) $ with weights $ w_n^{(B)}(0) $ defined by Eq.(\ref{init_weight_back_2}), we obtain the following formula for a large-$N$ approximation $ \bar{\Psi}_N(x,0) $ to the initial state $ \bar{\Psi}(x,0) $:
\beq
\label{bar_Psi_0_N}
\bar{\Psi}_N(x,0) = 
\frac{y^{1/4}}{2\pi}  \int_{\hat{y}}^{\infty} dz z^{\kappa -2} e^{-z/2} \left( \frac{\sin\left[ 2\sqrt{N}(\sqrt{z} - \sqrt{y} ) \right]}{
\sqrt{z} - \sqrt{y}} + \frac{\sin\left[ 2\sqrt{N}(\sqrt{z} + \sqrt{y} ) \right]}{
\sqrt{z} + \sqrt{y}} \right) 
\eeq
In the limit $ N \rightarrow \infty $, we use the relation
\beq
\label{limit_sin}
\lim_{N \rightarrow \infty} \frac{\sin\left[ 2\sqrt{N}(\sqrt{z} - \sqrt{y} ) \right]}{
\sqrt{z} - \sqrt{y}} = \pi \delta( \sqrt{z} - \sqrt{y} ) = 2 \pi \sqrt{y} \delta(z - y)
\eeq
to obtain 
\beq
\label{lim_bar_Psi_N}
\lim_{N \rightarrow \infty} 
\bar{\Psi}_N(x,0) = 
y^{\kappa - 5/4} e^{-y/2} \theta( y > \hat{y})
\eeq
The contribution of the second term in (\ref{bar_Psi_0_N}) vanishes in the limit $ N \rightarrow \infty $ as it is proportional to $ \delta(\sqrt{y} + \sqrt{\hat{y}}) $. Therefore we obtain 
\beq  
\label{lim_p_large_N}
\lim_{N \rightarrow \infty} 
p_N(x+ \hat{z},0) = \lim_{N \rightarrow \infty} e^{V(x+\hat{z})/\hbar} \bar{\Psi}_N(x,0) = \theta(y > \hat{y})
\eeq
consistently with the initial condition (\ref{init_cond_Psi}).  Therefore, we verified that the correct initial condition for $ p_N(x+ \hat{z},0) $ is recovered in the 
strict limit $ N \rightarrow \infty $. On the other hand, as in practice we have to work with a finite $ N $, it is important to estimate pre-asymptotic corrections to 
$ \bar{\Psi}_N(x,0) $ and $ p_N(x+ \hat{z},0) $ for large but finite values of $ N $. Calculation of these corrections is presented in Appendix B, and produces the following result (here $ \omega := 2 \sqrt{N} $):
 \beq
 \label{finite_bar_Psi_final} 
  \bar{\Psi}_{\omega}(x,0)  = 
 y^{\kappa - 5/4} e^{-y/2} \theta( y > \hat{y})  -
 \frac{2 y^{-1/4} \hat{y}^{\kappa - \frac{3}{2}} \sin \left( \omega \sqrt{\hat{y}} \right) }{\pi \omega}  
  \left[ 
  \cos \left( \omega \sqrt{y} \right) + 
  \sin \left( \omega \sqrt{y} \right) \right] + O \left(\frac{1}{\omega^{2}} \right) 
  \eeq

\begin{figure}[ht]
\begin{center}
\hspace*{-1cm} \includegraphics[
width=105mm,
height=65mm]{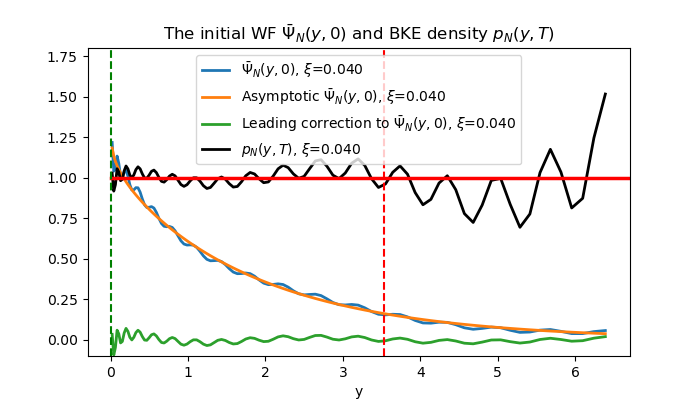}
\caption{  
Finite-$N$ approximation to the time-0 WF $\hat{\Psi}_N(y,0) $ and the terminal BKE tail probability $ p_N(y,T) $ computed according to (\ref{finite_bar_Psi_final}).
With $ N = 150 $, the finite-$N$ BKE density $ p_N(y,T) $ is an approximation to the true limit $ p_{\infty} (y,T) = \theta(y - \hat{y} )$ which is only recovered in the strict 
limit $
N \rightarrow \infty $.
 Increasing oscillations 
for large values of $ y $ indicate that corrections not included in  (\ref{finite_bar_Psi_final}) become important in this region. 
} 
\label{fig_BKE_density_and_Psi_bar_t0}
\end{center}
\end{figure}

The initial BK-SE WF $ \bar{\Psi}(y,0) $ and the terminal BKE density $ p_N(y,T) $ computed using Eq.(\ref{finite_bar_Psi_final}) are shown in 
Fig.~\ref{fig_BKE_density_and_Psi_bar_t0}.  A few comments can be made at this point. First, we observe that corrections to the asymptotic result given by the first term depend not just on the value of $N $ but rather on the combinations $ 2 \sqrt{N y} $ and $ 2 \sqrt{N \hat{y}} $. Second, 
corrections grow in the relative sense when $ y \rightarrow \infty $, showing that
Eq.(\ref{finite_bar_Psi_final_B}) is an asymptotic expansion 
valid only for sufficiently small values of $ y $ such that corrections are still smaller than the leading asymptotic term $ e^{-y/2} \theta( y > \hat{y}) $. Note that 
the leading correction in (\ref{finite_bar_Psi_final_B}) formally becomes a dominant term in the limit $ y \rightarrow \infty $, because it does not have an exponential term.
As the correct asymptotic behavior of the model at $ y \rightarrow \infty $ is actually given by the leading term  $ e^{-y/2} \theta( y > \hat{y}) $, it simply means that for larger values of $ y $ more and more terms in the eigenvalue decomposition become important. In the limit $ y \rightarrow \infty $, the whole infinite basis is required to produces the correct behavior.

Third, Fig.~\ref{fig_BKE_density_and_Psi_bar_t0} compares the finite-$N$ approximation to the BKE terminal density $ p_N(y,T) $ to its true asymptotic form $ 
\theta(y - \hat{y}) $. This analysis suggests that the BK-SE (\ref{SE_back_A}) may be more forgiving to a truncation of the infinite-dimensional basis than the related FP-SE (\ref{SE_A}). This can be based on the observation that at least with $ N = 150 $, our terminal BKE density $ p_N(y,T) $ looks more similar to 
its correct $ N \rightarrow \infty $ limit $ \theta(y > \hat{y}) $ than the time-0 FPE density $ f_N(y,0) $ looks similar to the delta-function $ \delta(y-y_0) $ in 
Fig.~\ref{fig_finite_N_FPE_t0}. Note that while the finite-$N$ density $ p_N(y,T) $ shows that corrections to the asymptotic expression become large for large values of $ y $, our approach is not expected to loose much accuracy due to the truncation as long as values of $ y $ of interest are still below this threshold.   

 Provided we use different initial weights given by (\ref{init_weight}) for the FP-SE (\ref{SE}) or (\ref{init_weights_back}) for the BK-SE (\ref{SE_back}), all subsequent steps are the same for both equations. 
 Substituting (\ref{eigen_val_exp}) into the SE (\ref{SE_A}), multiplying both sides 
 by $ \phi_m(x) $ and integrating over $ x $, we end up with a system of ordinary differential equations (ODEs) for coefficients $ w_n(\tau) $ which can be written 
 in a vector form as follows:
 \beq
 \label{ODE_w}
 - \frac{d {\bf w}(\tau)}{d \tau} = {\bf A}^{\pm} {\bf w}(\tau)
 \eeq 
 where we should $ {\bf A}_{+} $ or $ {\bf A}_{-} $ when $ g > 0 $ or $ g < 0$, respectively. The solution to this system of ODEs read
 \beq
 \label{ODE_solution}
 {\bf w}(\tau) = e^{- \tau {\bf A}^{\pm}} {\bf w}(0) = {\bf U}_{\pm} e^{- \tau \boldcalD_{\pm}} {\bf U}_{\pm}^T {\bf w}(0)
 \eeq
 where $ {\bf w}(0) $ is a vector of initial weights (\ref{init_weight}) or (\ref{init_weights_back}), and we used (\ref{A_eigenvalues}) in the second equation. The behavior of weights $ w_n(T) $ for different values of $ T = [5,10,15,20] $ is shown in Figs.~\ref{fig_weights_FPE_all_T} and \ref{fig_weights_BKE_all_T} for the FPE and BKE equations, respectively. As found in both cases, the dependence of weights $ w_n^{F)}(T) $ or $ w_n^{B)}(T) $ on $ n $ for long enough time horizons $ T $ is well approximated by the  formula $ w_n \simeq a e^{- b n} + c $ where parameters $ a,b,c $ depend on $ T $.

The time-dependent solution of the FP-SE (\ref{SE_A_0}) is now computed as follows:
\beq
\label{SE_res}
\Psi(x,t) = {\bf w}_{F}^{T}(0) e^{ - \frac{h}{2} t {\bf A}^{\pm}}  \boldphi(x)
\eeq
where we now denoted the FP-SE weights (\ref{init_weight}) as $ {\bf w}_{F}(0) $.
Now we can use Eqs.(\ref{Schrodinger_change_A}) and (\ref{exp_V_hbar}), we obtain the FPE density in the $ y $-space as follows:
 \beq
 \label{FPE_f_final_A}
f(y,t) = y^{-1} e^{ - V(y)/\hbar } \Psi(y,t) 
=  y^{\eta/\hbar-1} e^{ - y/2} \Psi(y,t)
 \eeq
Comparing (\ref{FPE_f_final_A}) with (\ref{SE_res}), we see that the 
FPE density $ f(z,t) $ has the same functional form as the quantum-mechanical WF (\ref{SE_res}), albeit with different coefficients, see 
Fig.~\ref{fig_FPE_and_BKE}. 

In its turn, the solution for BK-SE (\ref{SE_back_A}) is computed as follows:
\beq
\label{SE_back_res}
\bar{\Psi}(x,t) = {\bf w}_{B}^{T}(0) e^{ - \frac{h}{2} \left(T-t \right) {\bf A}^{\pm}}  \boldphi(x)
\eeq
where weights $ {\bf w}_{B}^{T}(0) $ are computed according to (\ref{init_weights_back}). 
Using Eqs.(\ref{Schrodinger_change_A}), we obtain the solution for the BKE: 
 \beq
 \label{BKE_f_final_A}
p(y,t) = e^{  V(y)/ \hbar } \Psi(y,t) 
=  y^{ -\eta/\hbar} e^{ y/2} \bar{\Psi}(y,t)
 \eeq 
Eqs.(\ref{SE_res}), (\ref{SE_back_res}), (\ref{FPE_f_final_A}), (\ref{BKE_f_final_A}) provide our final semi-analytical expressions for the quantum-mechanical WFs
$\Psi(x,t), \, \bar{\Psi}(x,t) $ and the FPE and BKE densities $ f(z,t), \, p(z,t) $, respectively. The behavior of the resulting densities is illustrated in Fig.~\ref{fig_FPE_and_BKE}.

\begin{figure}[ht]
\begin{center}
\hspace*{-2cm}  \includegraphics[
width=190mm,
height=65mm]{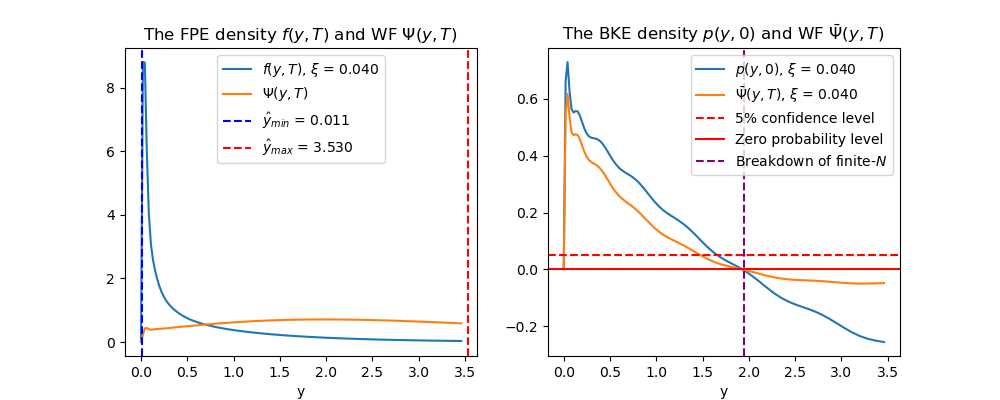}
\caption{On the left: The FPE density $ f(y, T) $ and quantum mechanical wave function $ \Psi(y,T) $ for $ T = 20 $. The FPE density $ f(y,T) $ has a sharp peak 
at a small value $ y= y_{\star} $, and monotonically decays for $ y > y_{\star} $. On the right:
the BKE tail probability $ p(y,0) $ and quantum mechanical wave function $ \bar{\Psi}(x,T) $  for $ T = 20 $. The BKE density has a maximum at a small value of $ y $, and monotonically decays to the right of this point. The vertical line denoted 'Breakdown of finite-$N$' shows the value of $ y $ where $ \bar{\Psi}(y,t) $ becomes slightly negative due to the basis truncation.}   
\label{fig_FPE_and_BKE}
\end{center}
\end{figure}

\def\thesection{B}	
\setcounter{equation}{0}
\def\theequation{\thesection.\arabic{equation}}

\section*{Appendix B: Asymptotic expansion of $ \bar{\Psi}_N(x,0) $}

In this appendix, we compute a leading pre-asymptotic correction to $ \bar{\Psi}_N(x,0) $ defined in Eq.(\ref{bar_Psi_0_N}). It is convenient at this stage to 
introduce the parameter $ \omega := 2 \sqrt{N} $, so that Eq.(\ref{bar_Psi_0_N}) takes the following form: 
\beq
\label{bar_Psi_0_N-B}
\bar{\Psi}_{\omega}(x,0) = 
\frac{y^{1/4}}{2\pi}  \int_{\hat{y}}^{\infty} dz z^{\kappa -2} e^{-z/2} \left( \frac{\sin\left[ \omega(\sqrt{z} - \sqrt{y} ) \right]}{
\sqrt{z} - \sqrt{y}} + \frac{\sin\left[ \omega(\sqrt{z} + \sqrt{y} ) \right]}{
\sqrt{z} + \sqrt{y}} \right) 
\eeq
To compute the pre-asymptotic correction, we use the identity
\beq
\label{integral_represent}
\bar{\Psi}_{\omega}(x,0) = \bar{\Psi}_{\infty}(x,0) - \int_{\omega}^{\infty} \frac{\partial \bar{\Psi}_{\omega'}(x,0)}{\partial \omega'} d \omega' := 
\bar{\Psi}_{\infty}(x,0) + \Delta_{\omega} \bar{\Psi}(x.0)
\eeq 
 where 
 \beq
 \label{psi_bar_infty}
\bar{\Psi}_{\infty}(x,0) := \lim_{\omega \rightarrow \infty} 
\bar{\Psi}_{\omega} (x,0) = 
y^{\kappa - 5/4} e^{-y/2} \theta( y > \hat{y})
\eeq
Therefore, we have to compute the asymptotics of the second term $ \Delta_{\omega} \bar{\Psi}(x.0) $ in (\ref{integral_represent}). Note that we have $ \lim_{\omega \rightarrow \infty} \Delta_{\omega} \bar{\Psi}(x,0) = 0 $ by construction.
To tackle the semi-infinite integration range $ [\omega, \infty] $ in (\ref{integral_represent}), we use the Bromwich integral representation of the Heaviside step-function\footnote{Here $ \varepsilon > 0 $ is arbitrary and can be set to zero at the end of calculations.}:
 \beq
 \label{theta_Brom}
 \theta(x) = \frac{1}{2 \pi i} \int_{\varepsilon - i \infty}^{\varepsilon + i \infty} dz \frac{e^{xz}}{z}
 \eeq
 along with Eqs.(3.462) and (9.246.3) from \cite{GR}. This produces the following expression for the derivative $ \partial {\bar \Psi}_N / \partial \omega $ in the limit 
 $ \omega \gg 1 $:
 \beq
 \label{der_Psi_bar}
 \frac{\partial \bar{\Psi}_N(x,0)}{\partial \omega} = - \frac{2 \Gamma(2 \kappa - 2)}{\pi} y^{\frac{1}{4}}
 \cos \left( \omega \sqrt{y} \right) \text{Re} \left[ \frac{e^{i \pi \kappa}}{2 \pi i} \int_{\varepsilon - i \infty}^{\varepsilon + i \infty} dz \frac{e^{- z \sqrt{\hat{y}}}}{z}
 \left[  \left( \omega - iz \right)^{2 - 2 \kappa} + \rho \left( \omega - iz \right)^{- 2 \kappa}  + \ldots \right] \right]
 \eeq
 where 
 $ \rho = (1-\kappa)(1- 2 \kappa) $, and ellipses stand for higher powers of $ 1/(\omega - iz) $. Plugging this back into (\ref{integral_represent}) and interchanging the order of integration with respect to $ \omega $ and $ z $, we obtain: 
 \beq
 \label{contour_aux}
 \Delta_{\omega} \bar{\Psi}(x.0) =  
 \frac{4 \Gamma(2 \kappa - 2)}{\pi} y^{1/4}
\text{Re} \left[ 
\frac{e^{i \pi \kappa}}{2 \pi i} \int_{\varepsilon - i \infty}^{\varepsilon + i \infty} dz \frac{e^{- z \sqrt{\hat{y}}}}{z}
 \left[ J \left(z,\omega, \nu_1 \right) + \rho J \left(z,\omega, \nu_2 \right)  + \ldots \right] \right]
 \eeq
 where $ \nu_1 = 2 - 2 \kappa $, $ \nu_2 = - 2 \kappa $, $ \kappa = q - \eta/\hbar $ and $ J(z,\omega, \nu) $ stands for the inner integral with respect to $ \omega' $: 
 \beq
 \label{J_w_z}
J(z,\omega) :=  \int_{\omega}^{\infty} d \omega' \cos \left(\omega' \sqrt{y} \right) (\omega' - i z)^{\nu} 
 \eeq
 The integral $ J(z, \omega, \nu) $ can be computed using Eqs.(3.761.2) and (3.761.7) from \cite{GR}:
 \beq
 \label{J_integral}
 J(z,\omega, \nu) = \frac{i y^{-\frac{\nu}{2} - \frac{1}{2}}}{2}  \left[ e^{-z \sqrt{y} - i \frac{\pi}{2} (\nu + 1)} \Gamma \left(\nu+1, i \sqrt{y}( \omega - i z) \right)
 - e^{z \sqrt{y} + i \frac{\pi}{2} (\nu + 1)} \Gamma \left(\nu+1, -i \sqrt{y}( \omega - i z) \right) \right]
 \eeq 
 While this expression gives an exact value of the integral for arbitrary $ \omega $, as we use asymptotic relations in other parts of our derivation, we can only retain the leading terms that survive in the limit $ \omega \rightarrow \infty $. This is achieved using the asymptotic expansion of the incomplete gamma function $ \Gamma(\nu, x ) $ given 
 by Eq.(8.357) in \cite{GR}:
 \beq
 \label{Gamma_asymptotics}
 \Gamma(\nu, x) = x^{\nu-1} e^{-x} \left[ 1 - (1-\nu) x^{-1} + \ldots \right], \; \; \; |x| \rightarrow \infty
 \eeq
  This gives the following asymptotic form for the integral $ J(z, \omega, \nu) $ valid in the limit $ \omega \gg 1 $:
  \beq
  \label{J_int_asymptotics}
  J(z, \omega, \nu) = - y^{-\frac{1}{2}} \left[ \left(\omega - iz \right)^{\nu} \sin \left( \omega \sqrt{y} \right) + \nu y^{-\frac{1}{2}} \left(\omega - iz \right)^{\nu-1} \cos \left( \omega \sqrt{y} \right) + O \left( \left( \omega - iz \right)^{\nu-2} \right) \right]
  \eeq  
Substituting this back into (\ref{contour_aux}), we obtain
\beq
\label{contour_delta_Psi}
 \Delta_{\omega} \bar{\Psi}(x.0) =   
 - \frac{2 \Gamma(2k - 2)}{\pi} y^{-1/4} \left[ \sin\left( \omega \sqrt{y} \right) 
 I(\hat{y}, \kappa) + 2 (1-\kappa) \cos \left( \omega \sqrt{y} \right)  I (\hat{y}, \kappa+\frac{1}{2} ) + \ldots \right]  
\eeq
where
\beq
\label{I_y_yhat_kappa}
I(\hat{y},\kappa) := \text{Re} \left[ \frac{e^{i\pi \kappa}}{2 \pi i} 
\int_{\varepsilon - i \infty}^{\varepsilon + i \infty} dz \frac{e^{- z \sqrt{\hat{y}}}}{z} \left( \omega - iz \right)^{2 - 2 \kappa} \right]
\eeq
This last integral can be transformed into a a closed contour integral by adding a semi-circle in the positive semi-plane at $ |z| \rightarrow \infty $.
The remaining contour integral can now be deformed to run from $ \infty $ to $ 0 $ above the branch cut on $ [0,\infty ] $, and then back from $ 0 $ to $ \infty $ under the branch cut. As the total value of the contour integral is zero, we obtain that the Bromwich integral arising in (\ref{contour_delta_Psi}) is given by the difference of two integrals $ I_{+} $ and $ I_{-} $ where the first integral runs from $ 0 $ to $ \infty $ with 
$ z = x $ with a real value of $ x $, while the second integral $ I_{-} $ runs from $ \infty $ to $ 0 $ with $ z = e^{2 i \pi} x $.
Due to the branch cut singularity, the two integrals do not cancel out. Instead, we find that $ I_{-} = - e^{4 i \pi \kappa} I_{+} $. 
This produces the following result for the integral (\ref{I_y_yhat_kappa}):
\beq
\label{I_y_yhat_res}
I(\hat{y},\kappa) = \frac{\sin(2 \pi \kappa)}{\pi} \left[ \frac{\Gamma(3 - 2 \kappa)}{\omega} \hat{y}^{\kappa - \frac{3}{2}} \sin \left( \omega \sqrt{\hat{y}} \right)
- \frac{\Gamma(4 - 2 \kappa)}{\omega^2} \hat{y}^{\kappa - 2} \cos \left( \omega \sqrt{\hat{y}} \right) \right]
\eeq
Plugging this into (\ref{contour_delta_Psi}) and using (\ref{integral_represent}), we finally obtain
\beq
 \label{finite_bar_Psi_final_B} 
 \bar{\Psi}_{\omega}(x,0)  = 
 y^{\kappa - 5/4} e^{-y/2} \theta( y > \hat{y})  -
 \frac{2 y^{-1/4} \hat{y}^{\kappa - \frac{3}{2}} \sin \left( \omega \sqrt{\hat{y}} \right) }{\pi \omega} 
  \left[ 
  \cos \left( \omega \sqrt{y} \right) + 
  \sin \left( \omega \sqrt{y} \right)  \right] + O \left(\omega^{-2} \right) 
  \eeq
We see that pre-asymptotic corrections die off very slowly $ \sim N^{-1/2} $ times oscillating $ \sin $ functions depending on $ 2 \sqrt{N y} $ and $ 2 \sqrt{N \hat{y}} $. Eq.(\ref{finite_bar_Psi_final_B}) is an asymptotic expansion 
valid only for sufficiently small values of $ y $ such that the correction is still smaller than the leading asymptotic term $ e^{-y/2} \theta( y > \hat{y}) $. Note that 
the correction in (\ref{finite_bar_Psi_final_B}) formally becomes a dominant term in the limit $ y \rightarrow \infty $, because it does not have a decaying exponential factor which appears in the first term.
As the correct asymptotic behavior of the model at $ y \rightarrow \infty $ is actually given by the leading term  $ e^{-y/2} \theta( y > \hat{y}) $, it simply means that for larger values of $ y $ more and more terms in the eigenvalue decomposition become important. In the limit $ y \rightarrow \infty $, the whole infinite basis is required to produce the correct behavior. Fortunately, the strict limit $ y \rightarrow \infty $ is mostly of academic interest, and in practice the actual values of $ y $ of interest for the BKE may lie well below the onset of such region.
 
 \def\thesection{C}	
\setcounter{equation}{0}
\def\theequation{\thesection.\arabic{equation}}

\section*{Appendix C: Probability current and constraints on parameter $ q $}

 In this appendix, we analyze probability conservation and implications of the basis truncation for the FPE and FP-SE, and show how this analysis is used in order to pick a proper value of our hyperparameter $ q $. 


The analysis below is based on another useful form for the FPE (\ref{FPE}) that expresses it in the form of a continuity equation
\beq
\label{FPE_prob_curr}
\frac{\partial f(x,t)}{\partial t} = \frac{\partial}{\partial x} J(x,t), \; \; \; J(x,t) := 
\frac{\partial V}{\partial x} f(x,t)  + \frac{\sigma^2}{2} \frac{\partial f(x,t)}{\partial x} 
\eeq
where $ J(x,t) $ is the probability current. Written in this form, the FPE shows that to conserve the total probability $ \int f(x,t) dx $, the probability current 
$ J(x,t) $ should vanish at $ x \rightarrow \pm \infty $.

It is instructive to start with verifying the conservation of probability for the expansion (\ref{eigen_val_exp}). According to (\ref{Schrodinger_change_A}), the 
normalization constraint $ \int f(x,t) dx = 1 $ for the FPE implies the following relation for the weights $ w_n^{(F)} $:
\beq
\label{normalization_for_Psi}
\int_{-\infty}^{\infty} f(x,t) dx = \int_{-\infty}^{\infty} e^{ - V(x)/\hbar } \Psi(x,t) dx = \sqrt{\Gamma(2q)} \sum_{n=0}^{\infty} w_N^{(F)}(\tau)
\int_{0}^{\infty}  y^{\eta/ \hbar - q - 1} \phi_0(y) \phi_n(y) dy  
\eeq
Here in the second equation we used the relation
\beq
\label{exp_V_hbar}
e^{- V(x)/\hbar} = y^{\eta/ \hbar} e^{- y/2} = \sqrt{\Gamma(2q)} y^{\eta/ \hbar - q} \phi_0(y)
\eeq
Comparing Eq.(\ref{normalization_for_Psi}) with the orthonormality condition $ \int dx \phi_n(x) \phi_m(x) = \int \frac{dy}{y} \phi_n(y) \phi_m(y) = \delta_{mn} $, it could be tempting to set $ q =  \eta/ \hbar $ here so that we would have
\beq
\label{would_have}
e^{- V(x)/ \hbar} \Psi(x,\tau) \sim \phi_0(y) \sum_{n=0}^{\infty} w_N^{(F)}(\tau) \phi_n(y)
\eeq
which would automatically preserve the total normalization due to orthogonality of the basis, as the right hand of this equation would then integrate to one provided the weights sum up to one.\footnote{Such expansions for the FPE density arise for problems with a stable ground state, see e.g. \cite{vanKampen}.}
However, as we will see shortly, we are {\it not} free to pick the value $ q = \eta/ \hbar $, as such choice produces a diverging probability current at infinity.
 
Back to Eq.(\ref{normalization_for_Psi}) and computing the integral in it, we obtain
\beq
\label{norm_check_2}
\int_{-\infty}^{\infty} f(x,t) dx =  
\frac{\Gamma(q + \eta/ \hbar)}{\Gamma(q - \eta/ \hbar)}
\sum_{n=0}^{\infty}  \frac{\Gamma(n+q - \eta/ \hbar)}{\sqrt{n! \Gamma(n + 2q)}} w_n(\tau)
\eeq
Note that this expression is linear, rather than quadratic in weights $ w_n(\tau) $, as would be the case in the conventional quantum mechanics \cite{Landau}.
A different normalization in our case is due to Eq.(\ref{Schrodinger_change_A}) which fixes the normalization rule for $ \Psi(x,\tau) $, and replaces here the 
standard quantum mechanical normalization rule $ \int \left| \Psi(x,\tau) \right|^2 dx = 1 $ which translates into the constraint $ \sum_{n} w_n^2(\tau) = 1 $ in 
terms of the eigenvalue expansion (\ref{eigen_val_exp}). 

The initial value of the total probability (\ref{norm_check_2}) can be computed explicitly using Eq.(\ref{init_weight}):
\bea
\label{normalization_for_Psi_0}
\int_{-\infty}^{\infty} f(x,0) dx  
&\hspace{-0.2cm} = \hspace{-0.2cm} &  y_0^{q - \eta/ \hbar} \frac{\Gamma(q + \eta/ \hbar)}{\Gamma(q- \eta/ \hbar)} \sum_{n=0}^{\infty} 
\frac{\Gamma(n + q  - \eta/ \hbar)}{\Gamma(2q+n)} L_n^{(2q-1)}(y_0) \nonumber \\
&\hspace{-0.2cm} = \hspace{-0.2cm} & y^{q-\eta/ \hbar}(x_0) y^{-q + \eta/ \hbar}(x_0) = 1 
\eea
Here in the second line we used the expansion of a power function $ f(x) = x^{-p} $ in Laguerre polynomials
\beq
\label{power_func_exp}
x^{-p} = \frac{\Gamma(\alpha+ 1 - p)}{\Gamma(p)} \sum_{n=0}^{\infty} \frac{\Gamma(n+p)}{\Gamma(n+ \alpha + 1)}  L_n^{(\alpha)}(x), \; \; \; (p> 0) 
\eeq
If we take $ p = q - \eta/\hbar $ in (\ref{power_func_exp}), we verify the last equation
in (\ref{normalization_for_Psi_0}). Importantly, as Eq.(\ref{power_func_exp}) requires that $ p > 0 $, we obtain that
\beq
\label{q_larger_eta2h}
q > \frac{\eta}{\hbar}
\eeq
This calculation shows that the FPE density is correctly normalized at time $ t = 0 $ due to Eq.(\ref{power_func_exp}), provided (\ref{q_larger_eta2h}) is satisfied.
While necessary, this condition is not sufficient on its own to provide a well-behaved model. In particular, we note that Eq.(\ref{q_larger_eta2h}) only gives a formal 
answer for an infinite sum. When the infinite sum in (\ref{power_func_exp}) is truncated to a finite sum as is done in practice, convergence properties
of the infinite series become important. 

For the series, (\ref{power_func_exp}), the convergence criterion amounts to the requirement that the integral $ \int_{0}^{\infty} x^{\alpha}e^{-x} x^{-2p} dx $ should stay finite. When we set here $ p = q - \eta/ \hbar $ and $ \alpha = 2 q - 1 $, this translates into the constraint $ \eta/\hbar > 0 $. However, this constraints does not always hold in our practical setting, where
typical values of $ \eta/ \hbar $ vary approximately in the range $ [-0.2, 0.2] $. When $ \eta/ \hbar < 0 $, the Taylor expansion (\ref{power_func_exp}) becomes 
an asymptotic series rather than a converging series. Unlike converging series, coefficients in asymptotic series may produce diverging patterns.    

 The asymptotic character of the Taylor expansion (\ref{power_func_exp}) with $ p = q - \eta/ \hbar $ when $ \eta/ \hbar < 0 $ can be seen as follows.   
The combinatorial factor $ \Gamma(n+ q - \eta/ \hbar)/  \Gamma(n + 2q) $ scales for $ n \rightarrow \infty $ as  
$ n^{-q - \eta/ \hbar} $, and the Laguerre polynomials $ L_N^{(2q-1)}(x) $ scale as $ n^{q-3/4} \sin \left( 2 \sqrt{nx} + \phi \right) $ where $ \phi $ is a phase parameter. Therefore, the coefficients in (\ref{normalization_for_Psi_0}) scale as
$ n^{-\eta/\hbar - 3/4} \sin \left( 2 \sqrt{n y_0} + \phi \right) $. Approximating the residual sum by an integral, we have an analytical expression for the residual of the infinite sum in (\ref{normalization_for_Psi_0}):
\bea
\label{P_more_than_N}
P_{>N} 
&\hspace{-0.2cm} := \hspace{-0.2cm} & y_0^{q -\eta/ \hbar} \frac{\Gamma(q + \eta/ \hbar)}{\Gamma(q- \eta/ \hbar)}  \sum_{n=N}^{\infty} 
\frac{\Gamma(n + q - \eta/ \hbar)}{\Gamma(2q+n)} L_n^{(2q-1)}(y_0) \nonumber \\
&\hspace{-0.2cm} = \hspace{-0.2cm} & 
y_0^{q -\eta/ \hbar} \frac{\Gamma(q + \eta/ \hbar)}{\Gamma(q- \eta/ \hbar)} \int_{N}^{\infty} x^{-\eta/\hbar -3/4} 
\sin \left( 2 \sqrt{y_0 x} + \phi \right) dx + \ldots
\eea
where ellipses stand for terms suppressed by powers of $ 1/N $.
This integral can be computed exactly in terms of incomplete gamma functions using Eqs.(3.761.2) and (3.761.7) in \cite{GR}, and it converges under the constraint
 \beq
 \label{tail_convergence_constraint}
 \eta/ \hbar > - 1/4  
 \eeq 
 In the limit $ \omega \rightarrow \infty $, the resulting expression considerably simplifies, and we obtain:
\beq
\label{P_large_N}
P_{>N} = 
\frac{\Gamma(q + \eta/ \hbar)}{\Gamma(q - \eta/ \hbar)} y_0^{q - 1/4} \left( y_0 N \right)^{- \eta/ \hbar -1/4} \cos \left( \sqrt{y_0 N} - \pi \left(q - 3/4 \right) \right) + 
O \left( \frac{1}{\omega^2} \right)  \; \; \;   \frac{\eta}{\hbar} > - \frac{1}{4}
\eeq
Fortunately the constraint $ \eta/ \hbar > - 1/4 $ is typically satisfied for cases of practical interest for our problem, where it ranges approximately in the interval $ [-0.2, 0.2] $. Also note that 
Eq.(\ref{P_large_N}) shows that for short times, the eigenvalue expansion is an asymptotic series rather than a converging series. The value of the residual sum $ P_{>N}  $ for a fixed $ N $ depends on the product $ N y_0 $ where $ y_0 $ is the initial value of the Morse variable, and 
Eq.(\ref{P_large_N}) implies that the value of $ N y_0  $ is such that the whole expression is between 0 and 1 - which would not be guaranteed for an arbitrary large value of $ N $ given the value of $ y $ due to the presence of the $ \cos $ function in (\ref{P_large_N}). 

When $ \eta/ \hbar > -1/4 $,
a finite-$ N $ approximation of the basis gives rise to strong oscillations up to moderately high values of $ N \sim 100-200 $, see Fig~\ref{fig_oscillating_norm}. These oscillations are due to in the oscillations of coefficients $ w_n(0) $ shown above in Fig.~\ref{fig_initial_weights}.  Furthermore, a finite-$ N $ approximation also introduces an artificial dependence of the norm of the FPE density on the initial value $ y_0 $ - which in fact exactly cancels out in (\ref{normalization_for_Psi_0}). 
 
\begin{figure}[ht]
\begin{center}
\hspace*{-1cm} \includegraphics[
width=85mm,
height=55mm]{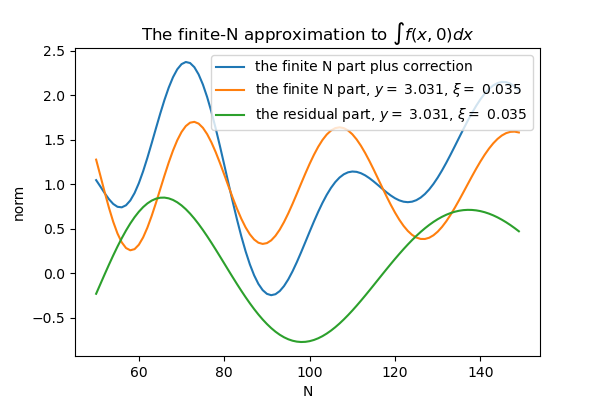}
\caption{Finite-$N$ approximation to the infinite series in Eq.(\ref{normalization_for_Psi_0}) as a function of $ N $. The blue line shows the result of a finite-$N$ truncation with a correction computed using Eq.(\ref{P_large_N}). In this figure, we chose to show the behavior for a large value of $ y_0 $, where corrections become of the same order as the leading term. For smaller values of  $ y_0 $, the leading correction is considerably smaller than the leading term, which justifies the use of asymptotic estimations in this case.
} 
\label{fig_oscillating_norm}
\end{center}
\end{figure}



We see that while the correct normalization of the FPE density at $ t = 0 $ is guaranteed in our theory, it critically depends on keeping the full infinite-dimensional basis. Related to this observation, 
and also implied by the principle of continuity, eigenvalue decomposition methods for the SE or FPE are known to perform poorly for short times due to 
a slow convergence, and thus are simply {\it not} a method of choice when the main focus of research is on a short-term behavior.
Fortunately, in our problem we are instead interested in a long-term behavior, where eigenvalue decomposition methods for PDEs usually work considerably better.
  

On the other hand, the problem still persists in practice.
A finite-$N$ truncation of the infinite basis needed for a numerical implementation creates initially unnormalized densities. 
As we have seen above, for long times of interest for this study $ T \simeq 15-20 $, weights decays with $ n $ very fast for both the forward and backward problems, as shown in Figs.~\ref{fig_weights_FPE_all_T} and \ref{fig_weights_BKE_all_T}. Nevertheless, the initial error in normalization due to the basis truncation 
propagates to time horizons of interest, at least with $ N = 150 $ (using $ N = 200 $ produced similar results), see Figs.~\ref{fig_norm_factor_for_FPE} and \ref{fig_norm_factors_various_T}. Norms of FPE densities for different time horizons 
show deviations from the unit normalization of the order of one percent for the range of parameters of interest. 

A properly normalized 
initial density can of course be obtained by dividing the unnormalized density by its total integral computed numerically. However, for our practical objective of comparing distributions  corresponding to different values of control parameters $ u_0 $ and $ \xi $, such approach can introduce some spurious behavior. In particular, 
our objective should be a monotonically increasing along both dimensions $ u_0 $ and $ \xi $. On the other hand, given that initial weights $ w_n(0) $ show a slow oscillatory decay, normalization of the initial FPE density $ f(x,0) $ by such procedure could introduce some non-monotonicity in the resulting behavior of our model.

\begin{figure}[ht]
\begin{center}
\hspace*{-1cm} \includegraphics[
width=175mm,
height=65mm]{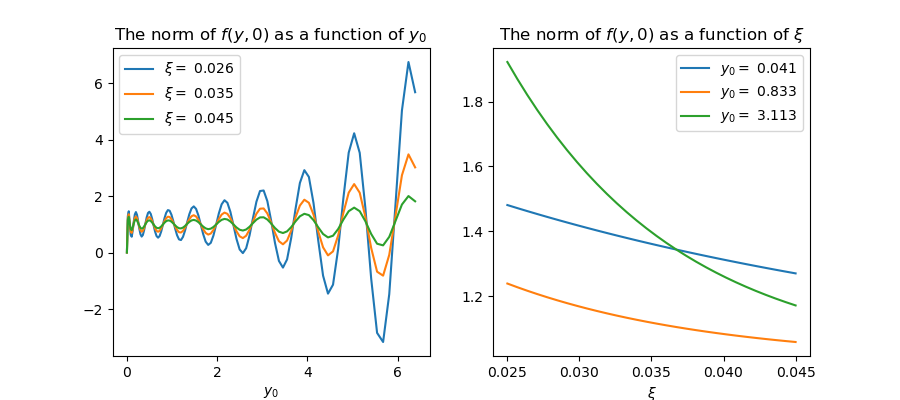}
\caption{The norm of the time-$0$ FPE density (\ref{FPE_f_final_A}) for $ T = 20 $, computed using $N=150$ basis functions. On the left: the norm as a function of initial contribution $ u_0 $ (shown here in units of the Morse variable y). One the right: the normalization factor as a function of contribution growth rate $ \xi $. Deviations of the norm from unity are due to the truncation of the infinite-dimensional basis to a finite basis of size $ N $. 
} 
\label{fig_norm_factor_for_FPE}
\end{center}
\end{figure}

\begin{figure}[ht]
\begin{center}
\hspace*{-1cm} \includegraphics[
width=175mm,
height=65mm]{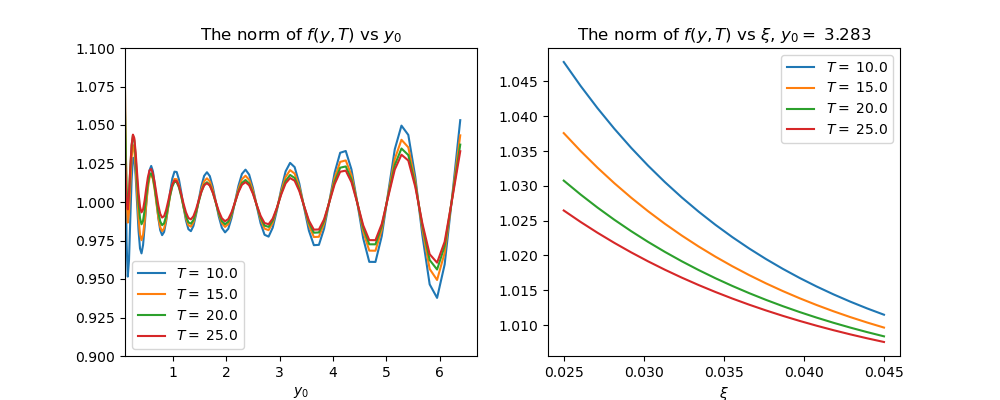}
\caption{The norm of the finite-$T$ FPE density (\ref{FPE_f_final_A}) for different values of $ T $, computed using $N=150$ basis functions. On the left: the norm as a function of initial contribution $ u_0 $ (shown here in units of the Morse variable $ y $), with a fixed value of $ \xi$. One the right: the norm as a function of contribution growth rate $ \xi $, with a fixed value of $ u_0 $. Deviations of the norm from unity are due to the truncation of the infinite-dimensional basis to a finite basis of size $ N $. 
} 
\label{fig_norm_factors_various_T}
\end{center}
\end{figure}
 
 To summarize to this point, pre-asymptotic corrections in expansions of FK-SE (\ref{SE}) and BK-SE (\ref{SE_back}) given, respectively, by Eqs.(\ref{P_large_N}) and (\ref{finite_bar_Psi_final_B}) show a similar behavior of both solutions as functions of $ 1/\sqrt{N} $, with oscillations and a slow power-like decay in the basis truncation number $ N $. When model parameters are taken to be in a range corresponding to real-world scenarios, we find that both the forward and backward models combined with the eigenvalue decomposition amount to asymptotic series in this regime. The similarity of the behavior of both WFs is of course not accidental as they are related by the Feynman-Kac relation (\ref{Psi_bar_Psi}). 
 
 In practice, when working with asymptotic series truncated at the $N$-th term, the important point to remember is that taking the largest possible value of $ N $ may not be the best idea. Instead, by analyzing relations such as  (\ref{P_large_N}) and (\ref{finite_bar_Psi_final_B}), and noticing that results depend on combinations of $ N $ with other inputs (the products $N y$ and $N y_0 $, in our case), a good value of $ N $ can be picked by the requirement that the pre-asymptotic correction should vanish, or that the partial sum in the expansion be stationary for an optimal value of $N $. Such tricks usually considerably improve models involving asymptotic expansions, see e.g. \cite{Migdal}. If our problem had fixed values of $ y $ or $ y_0 $, we could use (\ref{P_large_N}) and/or (\ref{finite_bar_Psi_final_B}) to pick a value of $ N $ such that the leading correction to the asymptotic formula would vanish. 
 Unfortunately, doing this could be problematic in the setting of the present work where the values of $ y $ and $ y_0 $ in Eqs.(\ref{P_large_N}) and (\ref{finite_bar_Psi_final_B}) are not fixed but instead serve as decision variables. Making $ N $ dependent on $ y $ or $ y_0 $ could introduce spurious non-monotonicity in the output BKE tail probability viewed as a function of control variables, and therefore we have not attempted such tricks for improving the behavior of asymptotic series.    
 
To verify that there is no non-vanishing probability flux at the spatial infinity, we compute the probability current $ J(z,t) $ as defined in Eq.(\ref{FPE_prob_curr}). Using Eq.(\ref{Schrodinger_change_A}), we obtain:
\beq
\label{J_z}
J(x+\hat{x},t) =  \frac{\hbar}{2} e^{- V(x+\hat{z})/\hbar} \left( \frac{1}{\hbar} \frac{\partial V}{\partial x} \Psi(x,t) + \frac{\partial \Psi}{\partial x} \right)
=  \frac{\hbar}{2} e^{- V(x+\hat{z})/\hbar} \mathcal{Q} \left(\eta/ \hbar \right) \Psi(x,t)
\eeq
To evaluate (\ref{J_z}), we use Eqs.(\ref{eigen_val_exp}), (\ref{raising_lowering}) and (\ref{exp_V_hbar}). This produces the following expression:
\beq
\label{J_z_2}
 J(y,t) =  \frac{\hbar}{2} 
  y^{\eta/ \hbar } e^{- y/2}
 \sum_{n=0}^{\infty} w_n(t) \left[ \mathcal{Q}(q+n)  -  \left( \frac{\eta}{\hbar}+ q +n \right) \right] \phi_n(y) 
 \eeq
 As the infinite sum in in this expression behaves in the limit $ y(x) \rightarrow 0 $ as $ y^{q} $, we obtain the $ \left. J(y,t) \right|_{y \rightarrow 0} \propto y^{ \eta/ \hbar+q} $. This shows that the requirement of vanishing probability current at $ y = 0 $ produces the following constraint for our hyperparameter $ q $:
 \beq
 \label{q_low_bound}
q > - \frac{\eta}{\hbar}
\eeq 
As the ratio $ \eta/\hbar $ varies in practical applications approximately in the range $ [-0.2, 0.2] $, the choice $ q = 5/4 $ made above in Eq.(\ref{q_final})
appears a satisfactory choice on all critically important model parameters (namely $ q $ and $ \eta $) given Eqs.(\ref{q_low_bound}) and (\ref{tail_convergence_constraint}). The behavior of probability current (\ref{J_z_2}) $ N = 150 $ basis functions is illustrated in Fig.~\ref{fig_prob_curr} for two times $ t = 0 $ and $ t = 20 $. 


\begin{figure}[ht]
\begin{center}
\hspace*{-1cm} \includegraphics[
width=175mm,
height=65mm]{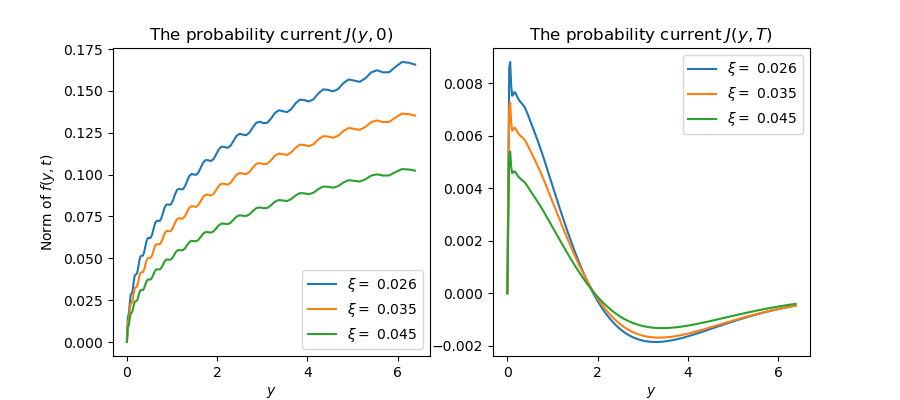}
\caption{The probability current $ J(y,t) $ according to Eq.(\ref{J_z_2})
using $ N = 150 $ basis functions. On the left: the initial probability current $ J(y,0) $. On the right: $ J(x,T) $  for $ T = 20 $.
} 
\label{fig_prob_curr}
\end{center}
\end{figure}

\end{document}